\documentclass[11pt]{amsart}
\usepackage{tikz-cd}
\usepackage{amsfonts}
\usepackage{mathrsfs}
\usepackage{amsxtra,amsthm,amsmath,amssymb}
\theoremstyle{plain}

\newcommand{\cleqn}{\setcounter{equation}{0}}
\newcommand{\clth}{\setcounter{theorem}{0}}
\newcommand {\sectionnew}[1]{\section{#1}\cleqn\clth}

\newtheorem{theorem}{Theorem}[section]
\newtheorem{lemma}[theorem]{Lemma}
\newtheorem{definition-theorem}[theorem]{Definition-Theorem}
\newtheorem{proposition}[theorem]{Proposition}
\newtheorem{corollary}[theorem]{Corollary}
\newtheorem{definition}[theorem]{Definition}
\newtheorem{example}[theorem]{Example}
\newtheorem{remark}[theorem]{Remark}
\newtheorem{notation}[theorem]{Notation}
\newtheorem{assumption}[theorem]{Assumption}
\newtheorem{lemma-definition}[theorem]{Lemma-Definition}
\newtheorem{lemma-notation}[theorem]{Lemma-Notation}
\newtheorem{question}[theorem]{Question}
\newtheorem{remark-definition}[theorem]{Remark-Definition}
\newcommand \bthm[1] { \begin{theorem}\label{t#1} }
\newcommand \ble[1] { \begin{lemma}\label{l#1} }

\newcommand \bpr[1] { \begin{proposition}\label{p#1} }
\newcommand \bco[1] { \begin{corollary}\label{c#1} }
\newcommand \bde[1] { \begin{definition}\label{d#1}\rm }
\newcommand \bex[1] { \begin{example}\label{e#1}\rm }
\newcommand \bre[1] { \begin{remark}\label{r#1}\rm }

\newcommand \bnota[1] {\begin{notation}\label{n#1}\rm }
\newcommand \bas[1] { \begin{assumption}\label{a#1}\rm }

\newcommand \bqu[1] { \begin{question}\label{q#1}\rm }

\newcommand {\ethm} { \end{theorem} }
\newcommand {\ele} { \end{lemma} }

\newcommand {\epr} { \end{proposition} }
\newcommand {\eco} { \end{corollary} }
\newcommand {\ede} { \end{definition} }
\newcommand {\eex} { \end{example} }
\newcommand {\ere} { \end{remark} }
\newcommand {\enota} { \end{notation} }
\newcommand {\eas} {\end{assumption}}

\newcommand {\equ} {\end{question}}
\newcommand \thmref[1]{Theorem \ref{t#1}}
\newcommand \leref[1]{Lemma \ref{l#1}}
\newcommand \prref[1]{Proposition \ref{p#1}}
\newcommand \coref[1]{Corollary \ref{c#1}}
\newcommand \deref[1]{Definition \ref{d#1}}
\newcommand \exref[1]{Example \ref{e#1}}
\newcommand \reref[1]{Remark \ref{r#1}}
\newcommand \lb[1]{\label{#1}}

\def \d {{\partial}}   

\def \Cset {{\mathbb C}}

\def \Pset {{\mathbb P}}

\def \cp   {{\Cset \Pset^1}}

\def \G  {{\mathcal{G}}}

\def \V {\mathcal V}

\def \al {\alpha}

\def \la {\lambda}

\def \ra  {\rightarrow}           

\def \lrw {\longrightarrow}

\def \la {\langle}
\def \ra {\rangle}

\def \g  {\mathfrak{g}}   
\def \h  {\mathfrak{h}}
\def \f  {\mathfrak{f}}
\def \n  {\mathfrak{n}}

\def \b  {\mathfrak{b}}

\def \sl {\mathfrak{sl}}

\def \d  {\mathfrak{d}}

\def \t  {\mathfrak{t}}

\def \c  {\mathfrak{c}}

\def \sC {{\scriptscriptstyle C}}

\def \sF {{\scriptscriptstyle F}}
\def \sG {{\scriptscriptstyle G}}

\def \sP {{\scriptscriptstyle P}}

\def \sX {{\scriptscriptstyle X}}
\def \sY {{\scriptscriptstyle Y}}
\def \sZ {{\scriptscriptstyle Z}}

\def \Ad { {\mathrm{Ad}} }

\DeclareMathOperator \ad { {\mathrm{ad}} }

\DeclareMathOperator \tr { {\mathrm{tr}} }

\DeclareMathOperator \Hom { {\mathrm{Hom}} }

\renewcommand \Im { {\mathrm{Im}} }

\newcommand{\beqa}{\begin{eqnarray*}}                     
\newcommand{\eeqa}{\end{eqnarray*}}
\def \hs {\hspace{.2in}}
\def \lara {\la \, , \, \ra}

\def \gog {\g \oplus \g}

\def \bs {{\bar{s}}}
\def \bu {{\bar{u}}}
\def \bv {{\bar{v}}}

\def \gst {\g_{\rm st}^*}
\def \gdia {\g_{{\rm diag}}}
\def \lam {\lambda}

\def \piP {\pi_{\scriptscriptstyle P}}
\def \piQ {\pi_{\scriptscriptstyle Q}}
\def \pist {\pi_{\rm st}}
\def \Pist {\Pi_{\rm st}}
\def \rst {r_{\rm st}}
\def \deltast {\delta_{{\rm st}}}

\def \sGs {{\scriptscriptstyle G}^*}
\def \piG {{\pi_{\scriptscriptstyle G}}}

\def \piGs {{\pi_{{\scriptscriptstyle G}^*}}}
\def \piX {{\pi_{\scriptscriptstyle X}}}
\def \piY {{\pi_{\scriptscriptstyle Y}}}

\def \piZ {\pi_{\sZ}}

\def \wF_mn {\wF_m \times \wF_n}
\def \wF_mnC {\wF_{m, n, \, \sC}}

\def \wF {\widetilde{F}}

\def \gotg {\g \otimes \g}
\def \Gvv {G^{v, v}}

\def \Guv {G^{u, v}}
\def \sSigma {{\scriptscriptstyle{\Sigma}}}
\def \ssG {{\scriptscriptstyle{\G}}}

\def \Suvt {{S}^{\bu, \bv}_{[t]}}

\begin{document}

\setlength{\baselineskip}{1.2\baselineskip}
\title[Double Bruhat cells and symplectic groupoids]{Double Bruhat cells and symplectic groupoids}
\author{Jiang-Hua Lu}
\address{
Department of Mathematics   \\
The University of Hong Kong \\
Pokfulam Road               \\
Hong Kong}
\email{jhlu@maths.hku.hk}
\author{Victor Mouquin}
\address{
Department of Mathematics   \\
University of Toronto \\
Toronto, Canada}               
\email{mouquinv@math.toronto.edu}
\date{}
\begin{abstract} 
Let $G$ be a connected complex semisimple Lie group, equipped with a standard multiplicative
Poisson structure $\pist$ determined by a pair of opposite Borel subgroups $(B, B_-)$.
We prove that for each $v$ in the Weyl group $W$ of $G$,
the double Bruhat cell $G^{v,v} = BvB \cap B_-vB_-$ in $G$, together with the Poisson structure $\pist$, is 
naturally a 
Poisson groupoid over the Bruhat cell $BvB/B$ in the flag variety $G/B$. Correspondingly,
every symplectic leaf of
$\pist$ in $\Gvv$ is a symplectic groupoid over $BvB/B$. For
$u, v \in W$, we show that the double Bruhat cell $(\Guv, \pist)$ has a naturally defined
left Poisson action by the Poisson groupoid 
$(G^{u, u}, \pist)$ and a right Poisson action by the Poisson groupoid $(\Gvv, \pist)$, and the two actions
commute. Restricting to symplectic leaves of $\pist$, one obtains commuting left and right Poisson actions 
on symplectic leaves in  $\Guv$ by symplectic leaves in $G^{u, u}$
and in $\Gvv$ as symplectic groupoids.
\end{abstract}
\maketitle


\sectionnew{Introduction and statements of results}\lb{intro}

\subsection{Introduction}
Let $G$ be a connected  complex semisimple Lie group, and let $(B, B_-)$ be a pair of
opposite Borel subgroups of $G$. It is well-known \cite{chari-pressley, etingof-schiffmann,
hodges, reshe-4, k-z:leaves} that the choice of $(B, B_-)$, together with that of a symmetric
non-degenerate invariant bilinear form on the Lie algebra of $G$,
determine a {\it standard multiplicative} Poisson structure $\pist$ on $G$
(see $\S$\ref{subsec-standard-G} for detail), and that the complex Poisson Lie group $(G, \pist)$
is the semi-classical limit of the quantized function algebra $\Cset_q[G]$ of $G$. 
The Poisson structure $\pist$ is invariant under left and right
translation by elements of the maximal torus $T = B \cap B_-$ of $G$, and it is well-known
\cite{hodges, reshe-4} that the {\it double Bruhat cells} 
$$
G^{u, v}=BuB \cap B_- v B_-,\hs \hs u, v \in W,
$$
where $W$ is the Weyl group of $(G, T)$, are precisely all the {\it $T$-leaves}
of $(G, \pist)$, i.e.,
submanifolds of $G$ of the form $\cup_{t \in T} \Sigma t$, where $\Sigma$ is a symplectic leaf of $\pist$ in $G$
(see \cite[$\S$2]{Lu-Mou:flags} on some basic facts of $\mathbb{T}$-leaves, where $\mathbb{T}$ is any torus). 
Double Bruhat cells have been
studied intensively and have served as motivating examples of the theories of total positivity and cluster algebras (see \cite{BFZ:clusterIII, Fomin-Zelevinsky:double, Goodearl-Yakimov:PNAS} and references therein).
When $G$ is simply connected, symplectic leaves of $\pist$ in each double Bruhat cell $G^{u,v}$ are
explicitly described in \cite{k-z:leaves}. 

The Poisson structure $\pist$ on $G$ projects to a well-defined Poisson structure $\pi_1$ on the
flag variety $G/B$, and each {\it Bruhat cell} $BvB/B \subset G/B$, where $v \in W$, is a Poisson
subvariety of $(G/B, \pi_1)$. In this paper, we show that for every $v \in W$ and any representative
$\bv$ of $v$ in the normalizer subgroup $N_G(T)$ of $T$ in $G$, the Poisson variety $(G^{v, v}, \pist)$
has a naturally defined groupoid structure over $BvB/B$, giving rise to a {\it Poisson groupoid} 
$(G^{\bv, \bv}, \pist)$ over the Poisson variety $(BvB/B, \pi_1)$.
The symplectic leaf $\Sigma^\bv$ of $\pist$ through $\bv$ is then shown to be a Lie sub-groupoid 
of $G^{\bv, \bv}$, becoming thus a symplectic groupoid over $(BvB/B, \pi_1)$.
The groupoid structure
on $G^{v, v}$ depends on the choice of $\bv \in N_G(T)$ (thus the notation $G^{\bv, \bv}$),
but different choices give isomorphic Poisson groupoids. For $u, v \in W$ and respective representatives
$\bu, \bv \in N_\sG(T)$, we show that the Bruhat cell $(\Guv, \pist)$ has a left Poisson action by
the Poisson groupoid $(G^{\bu, \bu}, \pist)$ and a right Poisson action by the Poisson groupoid
$(G^{\bv, \bv}, \pist)$, and the two actions commute. The two actions are then shown to restrict to 
commuting Poisson actions of the symplectic groupoids $(\Sigma^\bu, \pist)$ and $(\Sigma^\bv, \pist)$ on every 
symplectic leaf in $\Guv$.

\subsection{Statements of main results}\lb{subsec-Gvv-intro}
Let $v \in W$, and let $\bv$ be any representative of $v$ in 
$N_\sG(T)$. Let $C_\bv = N\bv \cap \bv N_-$, where $N$ and $N_-$ are respectively 
the uniradicals of $B$ and $B_-$. One then has the unique decompositions
$BvB = C_\bv B$ and $B_-vB_- = B_- C_\bv$ and the isomorphism
\[
C_\bv\, \stackrel{\sim}{\lrw} \,BvB/B, \;\;\;c \longmapsto c_\cdot B, \hs c \in C_\bv.
\]
Writing an element $g \in \Gvv$ uniquely as $g = c b = b_-c'$, where $b \in B$, $b_- \in B_-$, and $c, c' \in C_\bv$, the groupoid structure on $\Gvv$ over $BvB/B$ is defined as follows:
\begin{align*}
&\text{source map}: \; \theta_\bv(g) = g.B = c.B,   \\
&\text{target map}: \; \tau_\bv(g) = c'.B,    \\
&\text{inverse map}: \; \iota_\bv(g) = c'b^{-1} = b_-^{-1}c,   \\
& \text{identity bisection}: \; \epsilon_\bv(c_\cdot B) = c \in C_\bv \subset G^{v,v},\\
&\text{multiplication}:\; \mu_\bv(g, h) = c bb' = b_-b_-'c'', \; \mbox{if} \;
h =  c'b' = b_-'c'',\, \mbox{where}\; b' \in B, \, b'_- \in B_-, \, c'' \in C_\bv. 
\end{align*}
We will denote by  $G^{\bv, \bv} 
\rightrightarrows BvB/B$, or simply $G^{\bv, \bv}$, the 
double Bruhat cell $\Gvv$ with the groupoid structure thus defined. 
For another $u \in W$ and any  representative $\bu \in N_\sG(T)$ of $u$, define
\begin{align*}
\varpi:&\;\; \;\Guv \lrw BuB/B, \;\; \;\varpi(cb) = c_\cdot B, \hs \hs b \in B, \; c \in C_\bu,\\
\varpi': &\;\;\; \Guv  \lrw BvB/B,\;\;\; \varpi(b_-c') = c^\prime_\cdot B,
\hs \hs b_- \in B_-, \; c' \in C_\bv.
\end{align*}
The main results of the paper, \thmref{poiss-gpoid-Gvv}, \thmref{thm-Poi-bi}, and \thmref{symp-leaves}, 
can now be summarized as follows: let $u, v \in W$ and let $\bu$ and $\bv$ be any 
representatives of $u$ and $v$ in $N_\sG(T)$, respectively.

\vspace{.12in}
\noindent
{\bf Main Theorems.} {\it 1) The pair $(G^{\bv, \bv}, \, \pist)$ is a Poisson
groupoid over the Poisson manifold $(BvB/B, \pi_1)$, which, by restriction, also makes 
the symplectic leaf $\Sigma^\bv$ 
of $\pist$  through $\bv$ into a symplectic groupoid over $(BvB/B, \pi_1)$;

2) There is a natural left Poisson action of the Poisson groupoid $(G^{\bu, \bu}, \pist)$ on $(\Guv, \pist)$ with
moment map $\varpi$ and a natural right Poisson action of the Poisson groupoid $(G^{\bv, \bv}, \pist)$ on $(\Guv, \pist)$ with
moment map $\varpi'$. The two actions commute, and they restrict to Poisson actions of
the symplectic groupoids $(\Sigma^\bu, \pist)$ and $(\Sigma^\bv, \pist)$ on every symplectic leaf $\Sigma^{u, v}$ of
$\pist$ in $\Guv$.}

\vspace{.12in}
We remark that for any symplectic leaf $\Sigma^{u, v}$ in $\Guv$, the moment maps
\[
\varpi|_{\sSigma^{u, v}}:\;\;\; (\Sigma^{u, v}, \, \pist) \lrw (BuB/B, \, \pi_1)
\hs \mbox{and} \hs
\varpi'|_{\sSigma^{u, v}}: \;\;\; (\Sigma^{u, v}, \, \pist)  \lrw (BvB/B,\, -\pi_1)
\]
for the Poisson actions of
the symplectic groupoids $(\Sigma^\bu, \pist)$ and $(\Sigma^\bv, \pist)$ on  $(\Sigma^{u, v}, \, \pist)$
are symplectic realizations \cite{Wein:local, Xu:Morita} only 
in the sense that they are Poisson submersions but in general not
necessarily surjective (see \leref{le-Sigma-uv}, \leref{le-Guv-uv}, and \reref{re-Phi-uv-v}).

\vspace{.12in}
We in fact construct a Poisson groupoid $((G/B) \times B_-, \,\pi)$ over $(G/B, \pi_1)$, where the groupoid structure
is that of the action groupoid defined by the right action of $B_-$ on $G/B$ given by
\[
(g_\cdot B) \cdot b_- = (b_-^{-1}g)_\cdot B, \hs g \in G, \, \,b_- \in B_-,
\]
and the Poisson 
structure $\pi$ is a {\it mixed product Poisson structure} in the sense of \cite{Lu-Mou:mixed}, or, more precisely,
$\pi$ is the sum of the product Poisson structure $\pi_1 \times (\pist|_{\scriptscriptstyle{B_-}})$ on $(G/B) \times B_-$
and a certain mixed term determined by the action of $B$ on $G/B$
by left translation and by the action of $B_-$ on itself by left translation. 
For each $v \in W$ and a representative $\bv$ of $v$ in $N_\sG(T)$,
the Poisson groupoid $(G^{\bv, \bv}, \pist)$  
is then realized as a Poisson subgroupoid of the
Poisson groupoid $((G/B) \times B_-, \,\pi)$ over $(G/B, \,\pi_1)$ via a Poisson embedding
$I_\bv: (B_-vB_-, \pist) \to ((G/B) \times B_-, \,\pi)$
(see $\S$\ref{subsec-embedding} and $\S$\ref{subsec-Gvv-groupoid} for detail).
Using the embeddings $I_\bv$, we also interpret the Fomin-Zelevinsky twist map on double Bruhat cells
\cite{Fomin-Zelevinsky:double, k-z:leaves}
in terms of the inverse map of the groupoid $(G/B) \times B_-$ over $G/B$. See \reref{re-twist}.

The Poisson groupoid $((G/B) \times B_-, \pi)$ over $(G/B, \pi_1)$ is a special case of 
a general construction of action Poisson groupoids associated to quasitriangular $r$-matrices
(see $\S$\ref{subsec-transf-gpoid}). More precisely, given a Lie algebra $\g$, 
a quasitriangular $r$-matrix $r$ on $\g$, and a Lie algebra action of $\g$ on a manifold $Y$
such that the stabilizer subalgebra of $\g$ at each point of $Y$ is coisotropic with
respect to the symmetric part of $r$, Li-Bland and Meinrenken defined in \cite{Li-Mein} 
an {\it action Courant algebroid} over $Y$ with two transversal Dirac structures. In  
$\S$\ref{subsec-transf-gpoid}, we construct a pair of dual Poisson groupoids which integrate the two transversal 
Dirac structures in the sense that they have the two Dirac structures as their 
Lie bialgebroids (see \coref{co-F-pm} and \reref{re-F-pm} for detail). Applying the general construction to the semi-simple Lie algebra $\g$ and  
the standard quasitriangular $r$-matrix $r_{\rm st}$ on $\g$ (see $\S$\ref{subsec-standard-G}),
we obtain the action Poisson groupoid 
$((G/B) \times B_-, \pi)$ over $(G/B, \pi_1)$.

Symplectic groupoids and symplectic realizations are closely related to quantizations
of the Poisson manifolds \cite{Wein:local}. Relations between the symplectic groupoids of Bruhat cells described in
this paper and quantum Bruhat cells \cite{DC-K-P:quantum, L-M:quantum-Richardson, Lus:quantum, Yakimov:Spectra}
will be investigated in the future.

The paper is organized as follows. Some basic facts on Poisson Lie groups and Lie bialgebras are recalled in $\S$\ref{section-recalls}.
In $\S$\ref{sec-gpoids-algoids} we
construct a pair of dual action Poisson groupoids associated to quasitriangular $r$-matrices.
Some properties of the standard complex semisimple Poisson Lie groups are reviewed and proved in 
$\S$\ref{sec-review}. The main theorems of the paper are proved in 
$\S$\ref{sec-ss-PL} and $\S$\ref{sec-symplectic}, where we also generalize 
some results of \cite{k-z:leaves} on the symplectic leaves
of $\pist$ in the double Bruhat cells to the case when 
$G$ not necessarily simply connected.

\subsection{Acknowledgements}
This work was partially supported by the Research Grants Council of the Hong Kong SAR, China (GRF HKU 703712 and 17304415). 

\subsection{Notation}\lb{subsec-nota-intro}

Throughout this paper, vector spaces are understood to be real or complex. For a finite dimensional
vector space $V$, denote by $\lara$ the  
canonical pairing between $V$ and its dual space $V^*$. If $r = \sum_i x_i \otimes y_i \in V \otimes V$, let $r^{21} = \sum_i y_i \otimes x_i \in V \otimes V$ and let $r^\sharp: V^* \to V$ be the linear  map defined by 
$$
r^\sharp(\xi) = \sum_i \la \xi, x_i \ra y_i, \hs \xi \in V^*. 
$$

For a smooth (resp. complex) manifold $X$, the space of smooth (resp. holomorphic)
$k$-vector fields on $X$ will be denoted by $\V^k(X)$.
Let $X,Y$ be smooth or complex manifolds. For an integer $k \geq 1$ and 
$V_\sX \in \V^k(X)$ and $V_\sY \in \V^k(Y)$,
denote by $(V_\sX, 0)$ and $(0, V_\sY)$ the $k$-vector fields on $X \times Y$
whose values at $(x, y) \in X \times Y$ are respectively given by
$$
(V_\sX, \,0)(x, y) = i_y V_\sX(x) \hs \text{and} \hs (0, \,V_\sY)(x, y) = i_x V_\sY(y),
$$
where $i_y: X \to X \times Y, \,x' \mapsto (x', y)$ for $x' \in X$, and $i_x: Y \to X \times Y, \,y'\mapsto (x, y')$ for
$y' \in Y$.  We also denote $(V_\sX, 0) +(0, V_\sY)$ by $V_\sX \times V_\sY$.

Let $G$ be a Lie group with Lie algebra $\g$. A {\it left action} of $\g$ on a manifold $Y$
is a Lie algebra anti-homomorphism $\lam: \g \to \V^1(Y)$, while a {\it right action} of $\g$ on $Y$ is a Lie algebra homomorphism $\rho: \g \to \V^1(Y)$. If $\lam: G \times Y \to Y, (g, y) \mapsto gy$ is a
left action of $G$ on $Y$, one has the induced left action of $\g$ on $Y$, also denoted by $\lam$, given by 
$$
\lam:\; \g \longrightarrow \V^1(Y), \;\; \lam(x)_y = \frac{d}{dt}|_{t=0} \exp(tx)y, \hs x \in \g, \, y \in Y.
$$
Similarly, a right Lie group action $\rho: Y \times G \to Y, (y, g) \to yg$ induces a right Lie algebra action 
$$
\rho:\; \g \longrightarrow \V^1(Y), \;\;\rho(x)_y = \frac{d}{dt}|_{t=0} y\exp(tx), \hs x \in \g, \, y \in Y. 
$$
For $g \in G$, the left and right translation on $G$ by $g$, as well as their differentials, are
respectively denoted by $l_g$ and $r_g$.  If $k \geq 0$ is an integer and $x \in \g^{\otimes k}$, we denote by $x^L$ and $x^R$ the respective left- and right-invariant $k$-tensor fields on $G$ whose value at the identity element $e$ of $G$ is $x$. If $\xi \in \wedge^k \g^*$, we use similar notation for the left and right invariant $k$-forms with value $\xi$ at $e$.

Throughout the paper, if $(X, \pi)$ is a Poisson manifold and $X_1 \subset X$ a Poisson submanifold with respect to $\pi$, the
restriction of $\pi$ to $X_1$ will also be denoted by $\pi$ unless otherwise specified.

\sectionnew{Poisson Lie groups, $r$-matrices, and mixed product Poisson structures} \lb{section-recalls}

We recall from \cite{chari-pressley, etingof-schiffmann, Lu-Mou:mixed} some basic facts on Poisson Lie groups and
Lie bialgebras, and we refer to \cite[$\S$2]{Lu-Mou:mixed} in particular on certain conventions on constants and signs.

\subsection{Poisson Lie groups and Lie bialgebras} \lb{subsec-PL-Lie-bialg}
A {\it Lie bialgebra} is a pair $(\g, \delta_\g)$, where $\g$ is a (real or complex) finite dimensional Lie algebra, and $\delta_\g: \g \to \wedge^2 \g$ a linear map satisfying
$$
\delta_\g [x,y] = [x, \delta_\g (y)] + [\delta_\g(x), y], \hs x,y \in \g, 
$$
and such that the dual map $\delta_\g^*: \wedge^2 \g^* \to \g^*$ defines a Lie bracket on $\g^*$.
Given a Lie bialgebra $(\g, \delta_\g)$, the pair $(\g^*, \delta_{\g^*})$ is also a Lie bialgebra, 
where $\g^*$ is equipped with the Lie bracket dual to $\delta_\g$, and $\delta_{\g^*}: \g^* \to \wedge^2 \g^*$
is the dual map of the Lie bracket on $\g$.  One calls $(\g^*, \delta_{\g^*})$ the {\it dual Lie bialgebra} of $(\g, \delta_\g)$. If $(\g', \delta_{\g'})$ is any Lie bialgebra isomorphic to $(\g^*, \delta_{\g^*})$, we will
call $((\g, \delta_\g), (\g', \delta_{\g'}))$ a {\it pair of dual Lie bialgebras}.

Given a Lie bialgebra $(\g, \delta_\g)$, the vector space $\d = \g \oplus \g^*$ has a natural non-degenerate
symmetric bilinear form $\lara_\d$ defined by 
\begin{equation}\label{eq-lara-d}
\la x+ \xi, \; y + \eta \ra_\d = \la x, \,\eta\ra + \la y, \,\xi \ra, \hs x,\,y \in \g, \; \xi, \,\eta \in \g^*,
\end{equation}
and it is well-known that $\d$ has a unique Lie bracket $[\:, \:]$ such that both $\g$ and $\g^*$ are Lie
sub-algebras of $\d$ and such that $\lara_\d$ is ad-invariant. One calls $\d$ or $(\d, \lara_\d)$ the {\it double Lie algebra}
of $(\g, \delta_\g)$. Moreover, with $\delta_\d: \d \to \wedge^2 \d$ defined by 
$$
\delta_\d(x + \xi) = \delta_\g(x) - \delta_{\g^*}(\xi), \hs x \in \g, \; \xi \in \g^*,
$$
the pair $(\d, \delta_\d)$ is a Lie bialgebra, called 
 the {\it double Lie bialgebra} of $(\g, \delta_\g)$.

A {\it Poisson Lie group} is a pair $(G, \piG)$, where $G$ is a Lie group and $\piG$ a Poisson bivector field on $G$ 
that is {\it multiplicative} in the sense that the group multiplication $G \times G \to G$ is a Poisson map for the direct Poisson structure $\piG \times \piG$ on $G \times G$ and $\piG$ on $G$.  Given a Poisson Lie group $(G, \piG)$, the bivector field $\piG$ vanishes at the identity element $e$ of $G$, and the linearization $d_e \piG: \g \to \wedge^2 \g$ of $\piG$ at $e$, defined by $d_e \piG(x) = [\tilde{x}, \piG](e)$, where $\tilde{x}$ is any local vector field such that $\tilde{x}(e) = x$, is a Lie bialgebra structure on $\g$, and one calls $(\g, d_e \piG)$ the {\it Lie bialgebra of the Poisson Lie group} $(G, \piG)$. If $(G^*, \pi_{\sGs})$ is any Poisson Lie group whose Lie bialgebra is isomorphic to the dual Lie bialgebra of $(\g, d_e \piG)$, one says that $(G, \piG)$ and $(G^*, \pi_{\sGs})$ form a {\it pair of dual Poisson Lie groups}. 

Let $(G, \piG)$ be a Poisson Lie group with Lie bialgebra $(\g, \delta_\g)$, and let $\d$ be the double Lie algebra of $(\g, \delta_\g)$. Then $(G, \d)$ is a {\it Harish-Chandra pair} in the sense that the Lie algebra $\g$ of $G$ is a Lie subalgebra of $\d$,
and the Adjoint action $\Ad$ of $G$ on $\g$ extends to an action, still denoted by $\Ad$, of $G$ on $\d$ by 
Lie algebra automorphisms. Indeed, one has  \cite{dr:homog}
\begin{equation} \lb{lu-thm2.31}
\Ad_{g}\xi = r_{g^{-1}} \left(\pi_\sG^\#(g)(l_{g^{-1}}^* \xi)\right) + \Ad_{g^{-1}}^* \xi,\hs \xi \in \g^*,
\end{equation}
where  $\Ad_{g^{-1}}^*: \g^* \to \g^*$ is the dual map of $\Ad_{g^{-1}}: \g \to \g$ for $g \in G$. 

For $\xi \in \g^*$, the vector field ${\bf d}(\xi) = \pi_\sG^\#(\xi^R)$ on $G$ is called the
{\it dressing vector field} defined by $\xi$, where $\xi^R$ is the right invariant $1$-form on $G$
with value $\xi$ at $e$. By \eqref{lu-thm2.31}, one has
\begin{equation}\lb{eq-dressing}
{\bf d}(\xi)(g) = -l_g p_{\g}(\Ad_{g^{-1}} \xi), \hs \xi \in \g^*, \; g \in G,
\end{equation} 
where $p_\g: \d \to \g$ is the projection with respect to the decomposition $\d = \g + \g^*$. 



A {\it left Poisson action of a Poisson Lie group} $(G, \piG)$ on a Poisson manifold $(Y, \piY)$ is,
by definition, a left Lie group action $\lambda: G \times Y \to Y$ which is also a Poisson map
with respect to the product Poisson structure $\piG \times \piY$ on $G \times Y$ and the Poisson
structure $\piY$ on $Y$. Right Poisson actions of $(G, \piG)$ are similarly defined.
A {\it left Poisson action of a Lie bialgebra} $(\g, \delta_\g)$ on a Poisson manifold $(Y, \piY)$ is
a Lie algebra anti-homomorphism $\lambda: \g \to \V^1(Y)$ such that 
\[
[\lambda(x), \; \piY] = \lambda(\delta_\g(x)), \hs x \in \g,
\]
where $\lambda$ also denotes the linear map $\wedge^2 \g \to \V^2(Y)$ by 
$\lambda(x \wedge y) = \lambda(x) \wedge \lambda(y)$ for $x, y \in \g$. It is shown in 
\cite{Wein:remarks} that when a Poisson Lie group $(G, \piG)$ is connected, 
a Lie group action $\lambda: G \times Y \to Y$ of $G$ on a Poisson manifold $(Y, \piY)$
is a Poisson action of $(G, \piG)$ on $(Y, \piY)$ if
and only if the induced left Lie algebra action $\lambda: \g \to \V^1(Y)$
is a Poisson action of the Lie bialgebra $(\g, \delta_\g)$ of $(G, \piG)$ on
$(Y, \piY)$. Similar statement holds for right Poisson Lie group actions.

\subsection{Poisson structures defined by quasitriangular r-matrices} \lb{subsec-quasi-triangle}
Recall that a quasitriangular $r$-matrix on a Lie algebra $\g$ is
an element $r \in \g \otimes \g$ such that its symmetric part $\frac{1}{2} (r + r^{21})$ is invariant under the adjoint action of $\g$ on $\g \otimes \g$ and that $r$ satisfies the Classical Yang-Baxter Equation ${\rm CYB}(r) = 0$.
Given a quasitriangular $r$-matrix $r \in \gotg$, one has the Lie bialgebra $(\g, \delta_\g)$, where
$\delta_\g: \g \to \wedge^2 \g$ is defined by 
\begin{equation}\label{eq-delta-g-r}
\delta_\g(x) = \ad_x r, \hs x \in \g.
\end{equation}
A Lie bialgebra $(\g, \delta_\g)$ for which \eqref{eq-delta-g-r} holds for some quasitriangular $r$-matrix $r \in \gotg$ is said to be quasitriangular, and in such a case $r$ is called a {\it quasitriangular structure} of
$(\g, \delta_\g)$.

Let $r \in \gotg$ be a
quasitriangular $r$-matrix on a Lie algebra $\g$, and let $\sigma: \g \to \V^1(Y)$
be a right Lie algebra action of $\g$ on a manifold $Y$. 
If $r = \sum_i x_i \otimes x_i^\prime \in \gotg$, define 
\[
\sigma(r) = \sum_i \sigma(x_i) \otimes \sigma(x_i^\prime),
\]
which is a $2$-tensor field on $Y$. 
The following observation was made in \cite{Lu-Mou:mixed}.

\ble{le-sigma-pi}
If the $2$-tensor field $\sigma(r) \in \Gamma(TY \otimes TY)$ on $Y$ is skew-symmetric, then it is Poisson, and  
$\sigma$ is a (right) Poisson action of the Lie bialgebra $(\g, \delta_\g)$ on the
Poisson manifold $(Y, \sigma(r))$, where $\delta_\g$ is defined in \eqref{eq-delta-g-r}.
\ele

In the context of \leref{le-sigma-pi}, when $\sigma(r)$ is skew-symmetric,
the Poisson structure $\sigma(r)$ on $Y$ is said to be {\it defined} by the Lie algebra action $\sigma$ and the
$r$-matrix $r \in \gotg$.

\bre{re-s}
Let $s = \frac{1}{2} (r + r^{21})$ be the symmetric part of $r$. 
It is not hard to show (\cite[$\S$2.6]{Lu-Mou:mixed}) that $\sigma(r)$ is skew-symmetric, i.e., $\sigma(s) = 0$, if and only if the stabilizer subalgebra
of $\g$ at every $y \in Y$ is coisotropic with respect to $s$. Here
a subspace $\c$ of $\g$ is said to be {\it coisotropic with respect to $s$} if $s^\#(\c^0) \subset \c$, where 
$\c^0 = \{\xi \in \g^*: \la \xi, \c\ra = 0\} \subset \g^*$. 
\hfill $\diamond$
\ere

Let $(\g, \delta_\g)$ be a quasitriangular Lie bialgebra with quasitriangular $r$-matrix $r \in \g \otimes \g$.  Associated to $r$, one then \cite[$\S$2.3]{Lu-Mou:mixed} has the Lie subalgebras 
\begin{equation} \lb{defn f_pm}
\f_+ = \Im(r^\sharp) \hs  \text{and} \hs  \f_- = \Im((r^{21})^\sharp)
\end{equation}
of $\g$ and the Lie bialgebras $(\f_-, \delta_\g|_{\f_-})$ and $(\f_+, -\delta_\g|_{\f_+})$, which
are dual to each other under the pairing $\la \, \; \, \ra_{(\f_-, \f_+)}$ between $\f_-$ and $\f_+$ defined by
\begin{equation} \lb{f_pm-pairing}
\la (r^{21})^\sharp(\xi), \; r^\sharp(\eta)\ra_{(\f_-, \f_+)} = \la \xi, \; r^\sharp(\eta)\ra =
\la (r^{21})^\sharp(\xi), \; \eta \ra, \hs \xi,\; \eta \in \g^*. 
\end{equation}
If $(G, \piG)$ is a connected Poisson Lie group with Lie bialgebra $(\g, \delta_\g)$, and if
$F_+$ and $F_-$ are the connected Lie subgroups of $G$ with respective Lie algebras $\f_+$, $\f_-$,
then $F_+$ and $F_-$ are Poisson Lie subgroups of $(G, \piG)$. Moreover, denoting by the same symbol
the restrictions of $\piG$ to both $F_-$ and $F_+$,  $((F_-, \piG), (F_+, -\piG))$ is a pair of dual
Poisson Lie groups, with  $((\f_-, \delta_\g|_{\f_-}), (\f_+, -\delta_\g|_{\f_+}))$
as the corresponding pair of dual Lie bialgebras. 

\bex{exa-double-quasitriangle}
The double Lie bialgebra $(\d, \delta_\d)$ of any Lie bialgebra $(\g, \delta_\g)$ is quasitriangular, with
a quasitriangular structure defined by the quasitriangular $r$-matrix
\begin{equation} \lb{r-mat-double}
r_\d = \sum_{i=1}^n x_i \otimes \xi_i \in \d \otimes \d,
\end{equation} 
where $(x_i)_{i=1}^n$ is any basis of $\g$ and $(\xi)_{i=1}^n$ the dual basis of $\g^*$. 
In this example, the subalgebras $\f_+$ and $\f_-$ in \eqref{defn f_pm} are respectively $\g^*$ and $\g$. 
\hfill $\diamond$ 
\eex

\bre{re-df-g}
Let $(\g, \delta_\g)$ be a Lie bialgebra with a quasitriangular structure $r \in \g \otimes \g$.
Let $\d_{\f_-}$ be the double Lie algebra of $(\f_-, \delta_\g|_{\f_-})$, and let $r_{\d_{\f_-}} \in \d_{\f_-} \otimes \d_{\f_-}$ be defined as in \eqref{r-mat-double}.
Identifying $\f_-^* \cong \f_+$ via \eqref{f_pm-pairing}, the underlying vector space of $\d_{\f_-}$
is then $\f_- \oplus \f_+$, and the map
\begin{equation}\label{eq-q}
q: \;\; \d_{\f_-} \to \g, \hs q(x_-, x_+) = x_- + x_+, \hs x_- \in \f_-, \; x_+ \in \f_+, 
\end{equation}
is a Lie algebra homomorphism. Moreover (see \cite[Lecture 4]{etingof-schiffmann} and
\cite[$\S$2.3]{Lu-Mou:mixed}), $q(r_{\d_{\f_-}}) = r$. 
Thus if $(Y, \piY)$ is a Poisson manifold with a right Lie algebra action $\sigma: \g \to \V^1(Y)$ such that
$\piY$ is defined by $\sigma$ and $r$, i.e., $\piY = \sigma(r)$, then
$\piY$ is also defined by the Lie algebra $\d_{\f_-}$-action $\sigma \circ q: \d_{\f_-}
\to \V^1(Y)$ and the $r$-matrix $r_{\d_{\f_-}}$ on $\d_{\f_-}$.
\hfill $\diamond$
\ere

\subsection{Mixed product Poisson structures} \lb{subsec-mixed}
If $((\g, \delta_\g), (\g^*, \delta_{\g^*}))$ is a pair of dual  Lie bialgebras and if $(X, \piX)$ and $(Y, \piY)$ are Poisson manifolds with respective right and left Poisson actions
\[
\rho: \;\; \g^* \lrw \V^1(X) \hs
\mbox{and}\hs \lam: \;\; \g \lrw \V^1(Y)
\]
by Lie bialgebras, the bivector field $\piX \times_{(\rho, \lam)} \piY$ on the product manifold $X \times Y$ given by
\begin{equation}\label{eq-mixed}
\piX \times_{(\rho, \lam)} \piY = (\piX, 0) + (0, \piY) -\sum_{i=1}^n (\rho(\xi_i), \, 0) \wedge (0, \; \lam(x_i)),
\end{equation}
is a Poisson structure on $X \times Y$, called the {\it mixed product of $\piX$ and $\piY$ associated to $(\rho, \lam)$}, where $(x_i)_{i=1}^n$ is any basis for $\g$ and $(\xi_i)_{i=1}^n$ the dual basis for $\g^*$. 
We also refer to 
\[
-\sum_{i=1}^n (\rho(\xi_i), \, 0) \wedge (0, \; \lam(x_i)) \in \V^2(X \times Y)
\]
as {\it the  mixed term} of $\piX \times_{(\rho, \lam)} \piY$.
Mixed product Poisson structures of the form in \eqref{eq-mixed} are studied in \cite{Lu-Mou:mixed}.

\sectionnew{Action Poisson groupoids associated to quasitriangular $r$-matrices} \lb{sec-gpoids-algoids}

\subsection{Poisson groupoids} \lb{subsec-poiss-gpoids}
We recall from \cite{Mackenzie-Xu, Wein:Poi-goupoids, Xu:Poioid} some basic facts on Poisson groupoids.

Let $\G \rightrightarrows Y$ be a Lie groupoid, with $\theta, \tau: \G \to Y$ its source and target maps, $\iota: \G \to \G$ the groupoid inverse map, and  $\epsilon: Y \to \G$ the identity bisection. Let 
$$
\G_2 = \{ (a, b) \in \G \times \G : \tau(a) = \theta(b) \}
$$
be the submanifold of $\G \times \G$ of composable elements. A Poisson bivector field $\pi$ on $\G$ is
said to be {\it multiplicative} if  the graph of the groupoid multiplication 
$$
\{ (a, \,b, \, ab):\;  (a, b) \in \G_2 \} \subset  \G \times \G \times \G
$$
is a coisotropic submanifold of $\G \times \G \times \G$, where $\G \times \G \times \G$ is equipped with the Poisson structure $\pi \times \pi \times (-\pi)$. A {\it Poisson groupoid} is a pair $(\G \rightrightarrows Y, \pi)$, 
where $\G \rightrightarrows Y$ is a Lie groupoid and $\pi$ is a multiplicative Poisson structure on $\G$.
In such a case, $\iota(\pi) = -\pi$, and $\piY = \theta(\pi) = - \tau(\pi)$ is a Poisson structure
on $Y$, and one also says that $(\G \rightrightarrows Y, \pi)$ is a Poisson groupoid {\it over} $(Y, \piY)$. If in addition $\pi$ is non-degenerate,  one says that $(\G \rightrightarrows Y, \pi)$ is a {\it symplectic groupoid} over $(Y, \piY)$. 

Given a Lie groupoid $\G \rightrightarrows Y$, the left translation by $a \in \G$ is a smooth map
\[
l_a: \;\;\; \theta^{-1}(\tau(a)) \lrw \theta^{-1}(\theta(a)).
\]
Let $\ker \theta \to \G$ be the vector sub-bundle of the tangent bundle
of $\G$ whose fiber over
$a \in \G$ is the kernel of the differential of $\theta: \G \to Y$. A vector field $V$ 
on $\G$ is said to be {\it left invariant} if it is everywhere tangent to $\ker \theta$ and is invariant under the left translation
by every element in $\G$. The Lie algebroid of $\G \rightrightarrows Y$ is then the vector bundle $A = \epsilon^* \ker \theta$
over $Y$ with $\tau: A \to TY$ as the anchor map and with the Lie bracket on the space $\Gamma(A)$ of its sections defined by
\[
\stackrel{-\!\!\!-\!\!\!-\!\!\!-\!\!\!-\!\!\!-\!\!\!-\!\!\!\longrightarrow}{[s_1, s_2]} = 
[\stackrel{\longrightarrow}{s_1}, \; \stackrel{\longrightarrow}{s_2}],
\]
where for $s \in \Gamma(A)$, $\stackrel{\rightarrow}{s}$ is the unique left invariant vector field on $\G$
which coincides with $s$ on $\epsilon(Y) \cong Y$. As $T_{\epsilon(y)}\G =(\ker \theta)|_{\epsilon(y)} + T_{\epsilon(y)} \epsilon(Y)$
is a direct sum
for every $y \in Y$,  $A$ can be identified with the normal bundle of $\epsilon(Y)$ in $\G$.

If $(\G \rightrightarrows Y, \pi)$ is a Poisson groupoid, then 
the identity section $\epsilon(Y)$ is a coisotropic submanifold with respect to the Poisson structure $\pi$, and 
the dual vector bundle $A^*$ of $A$, identified with the co-normal bundle $N^*_{\epsilon(Y)}Y$
of $\epsilon(Y)$ in $\G$, is then a Lie sub-algebroid over $Y \cong \epsilon(Y) \hookrightarrow \G$
of the cotangent bundle Lie algebroid
$T_\pi^*\G$ over $\G$ defined by the Poisson structure $\pi$. The pair of Lie algebroids $(A, A^*)$ is then a
Lie bialgebroid \cite{Mackenzie-Xu} called the Lie bialgebroid of the Poisson groupoid 
$(\G \rightrightarrows Y, \pi)$. If $(\G' \rightrightarrows Y, \, \pi')$ is Poisson groupoid with Lie bialgebroid
$(A^*, A)$, one says that $((\G \rightrightarrows Y, \pi), \, (\G' \rightrightarrows Y, \, \pi'))$ is a {\it pair
of dual Poisson groupoids}.

Recall also that if $G$ is a Lie group and $\tau: Y \times G \to Y, (y, g) \to yg$, is a right Lie group action
of $G$ on a manifold $Y$, the product manifold $Y \times G$ then has the structure of an {\it action groupoid}, 
with $\tau: Y \times G \to Y$ as the target map, with 
\[
\theta(y, \, g) = y, \hs \hs y \in Y, \; g \in G, 
\]
as the source map, and with 
the groupoid multiplication, inverse map $\iota$, and the identity bisection $\epsilon$ respectively given by 
\begin{align*}
&(y_1, \, g_1)(y_2,\, g_2) = (y_1, \, g_1g_2), \hs \text{if} \hs y_1g_1 = y_2, \hs (y_1, \,g_1), \,(y_2, \,g_2) \in Y\times G,\\
&\iota(y, \, g) = (yg, \, g^{-1}), \hs \epsilon(y) = (y, \, e), \hs y \in Y, \; g \in G.
\end{align*}
Let ${\mathfrak{g}}$ be the Lie algebra of $G$.
Identifying $\epsilon^* \ker \theta$ with the trivial vector bundle
$A = Y \times \g$ over $Y$, the Lie algebroid of the action groupoid 
$Y \times G \rightrightarrows Y$ is then the {\it action Lie algebroid} $Y \times \g$ with anchor map, also denoted by $\tau$,
given by
\[
\tau: \;\; Y \times \g \lrw TY, \;\; \tau(y, x) = \frac{d}{dt}|_{t=0} y \exp(tx), \hs y \in Y, \, x \in \g,
\]
and the Lie bracket on its sections being the unique extending the Lie bracket on $\g$, identified with the space of
constant sections.  For $\varphi \in C^\infty(Y, \g) \cong \Gamma(Y \times \g)$, 
the left-invariant vector field $\stackrel{\rightarrow}{\varphi}$
on the action groupoid $Y \times G \rightrightarrows Y$ is then given by
\begin{equation}\label{eq-phi-arrow}
\stackrel{\rightarrow}{\varphi} (y, g) = (0, \; l_g \varphi(yg)), \hs y \in Y, \, g \in G.
\end{equation}
By an {\it action Poisson groupoid} we mean a Poisson groupoid whose underlying groupoid structure is that of an action 
groupoid.


\subsection{Action Poisson groupoids associated to quasitriangular $r$-matrices}\lb{subsec-transf-gpoid}
Let $(G, \piG)$ be a connected Poisson Lie group with Lie bialgebra $(\g, \delta_\g)$, and let $(\g^*, \delta_{\g^*})$ be the
dual Lie bialgebra of $(\g, \delta_\g)$. Let $(Y, \piY)$ be a Poisson manifold, and 
assume that $\rho: \g^* \to \V^1(Y)$ is a right Poisson action of the Lie bialgebra
$(\g^*, \delta_{\g^*})$ on $(Y, \piY)$.  One then has the mixed product Poisson structure $\pi$ on the product manifold 
$Y \times G$ given by
\begin{equation}\label{eq-pi-def}
\pi = \piY \times_{(\rho, \,\lam_{\sG})} \piG,
\end{equation}
where $\lam_\sG$ is the left Lie algebra action of $\g$ on $G$ generated by the left action of $G$ on itself by left multiplication, i.e., 
\[
\lam_\sG(x) = x^R, \hs x \in \g,
\]
where recall that for $x \in \g$, $x^R$ is the right invariant vector field on $G$ with value $x$ at the identity element $e$.
Assume that $G$ also acts on the right of $Y$ by
\[
\tau: \;\; Y \times G \lrw Y, \;\; (y, \, g) \longmapsto yg, \hs y \in Y,\, g \in G.
\]
Then $Y \times G$ has the corresponding structure of an
action groupoid over $Y$. We review in this section a necessary and sufficient condition for 
the pair $(Y\times G \rightrightarrows Y, \pi)$ to be a Poisson groupoid.

Let  $(\d, \delta_\d)$ be the double Lie bialgebra of $(\g, \delta_\g)$, where recall that $\d = \g \oplus \g^*$ as a vector space, and recall the quasitriangular $r$-matrix $r_\d$ on $\d$ given in \eqref{r-mat-double}. Let
\begin{equation}\lb{eq-sigma-tau-rho}
\sigma: \;\;\d \lrw \V^1(Y), \hs \sigma(x + \xi) = \tau(x) + \rho(\xi), \hs x \in \g, \; \xi \in \g^*,
\end{equation}
where $\tau$ also denotes the Lie algebra homomorphism $\g \to \V^1(Y)$ induced by the group action $\tau: Y \times G \to Y$
(see notation in $\S$\ref{subsec-nota-intro}). The following \thmref{transf-P-gpoid} was proved in \cite{lu:thesis}.

\bthm{transf-P-gpoid} \cite[Theorem 3.32]{lu:thesis}
The pair $(Y \times G \rightrightarrows Y, \;\pi)$ is a Poisson groupoid if
and only if  $\sigma: \d \to \V^1(Y)$ defined in \eqref{eq-sigma-tau-rho}
is a right Lie algebra action of $\d$ on $Y$ and 
$\piY= -\sigma(r_\d)$.
\ethm

As \cite{lu:thesis} is not published,
for the convenience of the reader,
we give an outline of the proof of \thmref{transf-P-gpoid} given in \cite{lu:thesis}.
We first prove a lemma which explains the main part of \thmref{transf-P-gpoid}.

For $\alpha \in \Omega^1(Y)$, let $X_\alpha = \pi^\#(\tau^*\alpha) \in \V^1(Y \times G)$.
By \cite[Proposition 2.7]{Xu:Poioid}, if $(Y \times G \rightrightarrows Y, \;\pi)$ is a Poisson groupoid,  
$X_\alpha$ is necessarily a left invariant vector field on $Y \times G$ for every $\alpha \in \Omega^1(Y)$,
i.e., $\theta(X_\alpha) = 0$ and $X_\alpha(ab) = l_a X_\alpha(b)$ for any composable pair $(a, b)$ in $Y \times G$. 

\ble{le-X-alpha-0} 
1) One has $\theta(X_\alpha) = 0$ for all $\alpha \in \Omega^1(Y)$ if and only if $\piY = -\sigma(r_\d)$;

2) Assume that $\piY = -\sigma(r_\d)$. Then $X_\alpha$ is left invariant for all $\alpha \in \Omega^1(Y)$ if and only if  $\sigma: \d \to \V^1(Y)$ is a right Lie algebra action.  In such a case, for
$\alpha \in \Omega^1(Y)$, one has
$$
X_{\alpha} = \stackrel{\longrightarrow}{\phi_\alpha},
$$
where $\phi_\alpha \in C^\infty(Y, \g)$ is given by $\phi_\alpha(y) = -\rho_y^*(\alpha(y))$, with
$\rho_y: \g^* \to T_yY$ given by $\rho_y(\xi) = \rho(\xi)(y)$ for 
$y \in Y$ and $\xi \in \g^*$. 
\ele

\begin{proof}
For  $g \in G$ and $y \in Y$, let
\[
\tau_g: \;\; Y \lrw Y, \;\; y' \longmapsto y'g, \hs y' \in Y, \hs \mbox{and} \hs 
\tau_y:\;\; G \lrw Y, \;\; g' \longmapsto yg', \hs g' \in G.
\]
Let $p_1: Y \times G \to Y$ and $p_2: Y \times G \to G$ be respectively the projections to the first and
the second factors. Let $\alpha \in \Omega^1(Y)$ and let $y \in Y$ and $g \in G$. Then  
\[
(\tau^*\alpha)(y, \,g) =p_1^* \tau_g^* \alpha(yg)  + p_2^* l_{g^{-1}}^* \tau_{yg}^* (\alpha(yg)) \in T_{(y, g)}^*(Y \times G).
\]
Using the definition of $\pi$, one has
\begin{equation}\lb{eq-pi-tau}
X_\alpha(y, \, g) = (\pi_\sY^\#(y)(\tau_g^*\alpha(yg))+\rho_y(\tau_y^*\tau_g^*\al(yg)), \;\,
\pi_\sG^\#(g) (l_{g^{-1}}^* \tau_{yg}^* \al(yg)) - r_g \rho_y^* \tau_g^* \al(yg)).
\end{equation}

1) Let $\{x_i\}_{i=1}^n$ be any basis of $\g$ and $\{\xi_i\}_{i=1}^n$ the dual basis of $\g^*$, so that
$r_\d = \sum_{i=1}^n x_i \otimes \xi_i \in \d \otimes \d$. Then
$\piY = -\sigma(r_\d)$ if and only if $\piY =- \sum_{i=1}^n \tau(x_i) \otimes \rho(\xi_i)$,
which is equivalent to
\[
\pi_\sY^\#(y) (\alpha_y) = -\rho_y (\tau_y^* \alpha_y), \hs \hs y \in Y, \; \alpha_y \in T_y^*Y.
\]
It is now clear from  \eqref{eq-pi-tau} that 1) holds.

2) Assume now that $\piY = -\sigma(r_\d)$. By \eqref{eq-pi-tau},  $X_\alpha$ is left invariant if and only if 
\begin{equation} \lb{iif-left-inv}
-l_g \rho_{yg}^*\al(yg) =  \pi_\sG^\#(g) (l_{g^{-1}}^* \tau_{yg}^* \al(yg)) - r_g \rho_y^* \tau_g^* \al(yg), \hs 
(y, \, g) \in Y \times G. 
\end{equation}
Pairing both sides of \eqref{iif-left-inv} with $l_{g^{-1}}^*\xi \in T_g^*G$, where $\xi \in \g^*$, and using \eqref{lu-thm2.31}, one can rewrite \eqref{iif-left-inv} as 
$$
0 = \la \alpha(yg), \;\; \rho_{yg}(\xi) - \tau_g \tau_y(p_\g(\Ad_g \xi)) - \tau_g \rho_y(p_{\g^*}(\Ad_g \xi)) \ra,
\hs y \in Y,\, g \in G,
$$
where recall that $p_\g: \d \to \g$ and $p_{\g^*}: \d \to \g^*$ are the projections with respect to the decomposition
$\d = \g + \g^*$. Therefore $X_\alpha$ is left-invariant for all $\alpha \in \Omega^1(Y)$ if and only if 
\begin{equation} \lb{iif-left-inv-2}
\tau_{g^{-1}}(\rho(\xi)) = \sigma(\Ad_g \xi) \in \V^1(Y), \hs g \in G, \; \xi \in \g^*. 
\end{equation}
Assuming \eqref{iif-left-inv-2} and differentiating $g \in G$ in the direction of $x\in \g$ gives
\begin{equation}\label{eq-sigma-anti}
[\tau(x), \, \rho(\xi)] = \sigma([x, \, \xi]), \hs x \in \g, \, \xi \in \g^*,
\end{equation}
so $\sigma: 
\d \to \V^1(Y)$ is a Lie algebra homomorphism.
Conversely, assume that $\sigma$ is a Lie algebra homomorphism. 
The infinitesimal $\g$-invariance of $\sigma$ in \eqref{eq-sigma-anti} and the connectedness of $G$ imply the
$G$-equivariance of the $\sigma$, namely \eqref{iif-left-inv-2}. It is also clear 
from \eqref{eq-pi-tau} that in such a case, $X_\alpha =
\stackrel{\longrightarrow}{\varphi_\alpha}$ with $\varphi_\alpha$ as described.
\end{proof}

\noindent
{\it Proof of \thmref{transf-P-gpoid}}
Assuming that $\sigma: \d \to \V^1(Y)$ is a Lie algebra homomorphism and that
$\piY = -\sigma(r_\d)$, we now show that $(Y \times G \rightrightarrows Y, \pi)$ is a Poisson groupoid,
the other direction of \thmref{transf-P-gpoid} having been proved in \leref{le-X-alpha-0}.


Note first that $\piY = -\sigma(r_\d)$ implies that $\sigma$ is a right Poisson action of 
the double Lie bialgebra
$(\d, -\delta_\d)$ on $(Y, \piY)$, so $\tau$ is a right Poisson action of
the Poisson Lie group $(G, -\piG)$ on $(Y, \piY)$. It follows by an easy calculation that
the target map $\tau: (\G, \; \pi) \to (Y, \, \piY)$ 
is anti-Poisson, where $\G = Y \times G$. As the source map $\theta: (\G, \; \pi) \to (Y, \, \piY)$ is Poisson, the submanifold
\[
\G_2 = \{(a, \, b) \in \G \times \G: \tau(a) = \theta(a)\}\subset \G \times \G
\]
of composable pairs is coisotropic with respect to the product Poisson structure $\pi \times \pi$ on 
$\G \times \G$. Note that the graph 
$\{(a, \,b, \, ab):  (a, b) \in \G_2 \} \subset \G \times \G \times \G$ of the groupoid multiplication is the graph of the map
$\mu|_{\scriptscriptstyle{\mathcal{G}_2}}: \G_2 \to \G$, where
\[
\mu: \;\; \G \times \G \lrw \G: \;\; (y_1, \, g_1, \, y_2, \, g_2) \longmapsto (y_1, \, g_1g_2), 
\hs y_i \in Y, \, g_i \in G.
\]
Note also that the map $\mu$ is Poisson with respect to the product Poisson structures $\pi \times \pi$ on $\G \times \G$ and $\pi$ on $\G$. Indeed, $\mu = \nu \circ ({\rm Id}_{\scriptscriptstyle{\mathcal{G}}} \times p)$, where
the projection $p: (Y \times G, \, \pi) \to (G, \,\piG)$ to the second factor
is Poisson, and the map
\[
\nu: \; (Y \times G, \; \pi) \times (G, \, \piG) \lrw (Y \times G, \, \pi), \;\; (y, \, g, \, g_1) \longmapsto (y, \, gg_1),
\hs y \in Y,\, g, \, g_1 \in G. 
\]
is Poisson. It is a general fact, the  proof of which is straightforward (see \cite[Lemma 3.33]{lu:thesis}),  
that for a Poisson map  $\Phi: (P, \piP) \to (Q, \piQ)$ and 
a coisotropic submanifold $P_1 \subset (P, \piP)$, the graph $\{(p,\, \Phi(p)): p \in P_1\}$ 
of $\Phi|_{P_1}: P_1 \to Q$ is 
coisotropic in $(P \times Q, \, \piP \times (-\piQ))$ if and only if 
\[
\pi_\sP^\#(N_{\sP_1}^*P) \subset \ker \Phi,
\]
where $N_{\sP_1}^*P \subset T^*P|_{\sP_1}$ is the co-normal space of $P_1$ in $P$, and the sub-bundle 
$\pi_\sP^\#(N_{\sP_1}^*P)$ of $TP_1$ is called the {\it characteristic distribution} of the coisotropic submanifold
$P_1$ in $P$. Using \leref{le-X-alpha-0}, a direct calculation shows that the characteristic distribution
of $\G_2$ in $\G \times \G$ at the point $(y_1, \, g_1, \, y_2, \, g_2) \in \G_2$ is given by 
\[
\{(0, \; -l_{g_1} x, \; 
-\tau_{y_2} x, \; r_{g_2} x): \; x \in \rho_{y_2}^* T_{y_2}^* Y \subset \g\},
\]
which is easily seen to be contained in the kernel of the differential of $\mu$ at $(y_1, \, g_1, \, y_2, \, g_2) \in \G_2$.
Thus the graph $\{(a, \,b, \, ab):  (a, b) \in \G_2 \}$ is a coisotropic submanifold of
$\G \times \G \times \G$ with respect to the Poisson structure $\pi \times \pi \times (-\pi)$, and hence
$(\G \rightrightarrows Y, \pi)$ is a Poisson groupoid. 
This finishes the proof of \thmref{transf-P-gpoid}.

\bre{re-pair} In the context of \thmref{transf-P-gpoid}, it is easy to see that
the Lie algebroid structure induced by $\pi$ on the co-normal bundle of $\epsilon(Y)$
in $Y \times G$, identified with the trivial vector bundle $Y \times \g^*$ over $Y$, 
is that of the action Lie algebroid defined by the
right action $\rho$ of $\g^*$ on $Y$. Thus the Lie bialgebroid of the Poisson groupoid 
$(Y \times G \rightrightarrows Y, \pi)$ is the pair 
\[
(A = Y \times \g, \; \;A^* = Y \times \g^*)
\]
of action Lie algebroids.  Their double, as a Courant Lie bialgebroid
\cite{Liu-Wei-Xu}, 
is  the {\it action Courant algebroid} $Y \times \d$ over $Y$ defined by $\sigma$ 
that has been studied in \cite{Li-Mein}. 
\hfill $\diamond$
\ere

Let $((G, \piG), \, (G^*, \piGs))$ be a pair of Poisson Lie groups, with 
$((\g, \delta_\g), (\g^*, \delta_{\g^*}))$ the corresponding pair
of dual Lie bialgebras, and let $(\d, \lara_\d)$ be their
double Lie algebra. Let again $r_\d = \sum_{i=1}^n x_i \otimes \xi_i$ be the quasitriangular $r$-matrix on $\d$, where
$\{x_i\}_{i=1}^n$ is any basis of $\g$ and $\{\xi_i\}_{i=1}^n$ the dual basis of $\g^*$.
Assume that $\sigma: \d \to \V^1(Y)$ is a right Lie algebra action of $\d$ on a manifold $Y$ such that 
the stabilizer subalgebra $\d_y$ of $\d$ at every $y \in Y$ is coisotropic with respect to $\lara_\d$, which, by
\reref{re-s}, is equivalent to $\sigma(r_\d)$ being a Poisson structure on $Y$.
\bco{co-GG}
1) Assume that $\sigma|_\g: \g \to \V^1(Y)$ integrates to a Lie group action $Y \times G \to Y$. Then 
one has the action Poisson groupoid $(Y \times G \rightrightarrows Y, \, \pi_{\sY \times \sG})$ over $(Y, -\sigma(r_\d))$, where
$Y \times G \rightrightarrows Y$ is the action groupoid over $Y$ defined by the group action of $G$ on $Y$, and
$\pi_{\sY \times \sG}$ is the mixed product Poisson structure on $Y \times G$ given by
\[
\pi_{\sY \times \sG} = (-\sigma(r_\d), \, 0) + (0, \, \piG) - \sum_{i=1}^n \sigma(\xi_i), \, 0) \wedge (0, \, x_i^R);
\]

2) Assume that $\sigma|_{\g^*}: \g^* \to \V^1(Y)$ integrates to a Lie group action $Y \times G^* \to Y$. Then 
one has the action Poisson groupoid $(Y \times G^* \rightrightarrows Y, \, \pi_{\sY \times \sG^*})$
over $(Y, \sigma(r_\d))$, where
$Y \times G^* \rightrightarrows Y$ is the action groupoid over $Y$ defined by the group action of $G^*$ on $Y$,  and
$\pi_{\sY \times \sG^*}$ is the mixed product Poisson structure on $Y \times G$ given by
\[
\pi_{\sY \times \sG^*} = (\sigma(r_\d), \, 0) + (0, \, \piGs) - \sum_{i=1}^n \sigma(x_i), \, 0) \wedge (0, \, \xi_i^R);
\]

3) When the assumptions in both 1) and 2) hold, the two action Poisson groupoids in 1) and 2) form a dual pair of Poisson
groupoids.
\eco

\begin{proof}
By  \leref{le-sigma-pi}, $\sigma$ is a right Poisson action of the Lie bialgebra $(\d, \delta_\d)$ on $(Y, \sigma(r_\d))$,
where recall that $\delta_\d(v) = \ad_v r_\d$ for $v \in \d$. As $\delta_\g = \delta_\d|_\g$ and $\delta_{\g^*} = 
-\delta_\d|_{\g^*}$, $\sigma|_{\g^*}$ is a right Poisson action of the Lie bialgebra $(\g, \delta_{\g^*})$ on 
$(Y, -\sigma(r_\d))$, and $\sigma|_\g$ is a right Poisson action of the Lie bialgebra $(\g, \delta_\g)$ on $(Y, \sigma(r_\d))$.
Now 1) and 2) of \coref{co-GG} follows from 
\thmref{transf-P-gpoid} applied to the Poisson Lie groups $(G, \piG)$ and $(G^*, \piGs)$ respectively, and
3) follows from \reref{re-pair}.

\end{proof}

\bre{re-charish-chandra}
When $\sigma|_\g: \g \to \V^1(Y)$ integrates to a Lie group action $\tau: Y \times G \to Y$, the pair $(\tau, \sigma)$
can be thought of as a (right) {\it action of the Harish-Chandra pair}
$(G, \d)$ (see $\S$\ref{subsec-PL-Lie-bialg}) on the manifold $Y$ in the sense that 
$\tau$ is 
a right action of the Lie group $G$ on $Y$ and $\sigma$ is a right action of the Lie algebra $\d$ on $Y$
such that $\sigma|_\g$ coincides with the action of $\g$ on $Y$ induced by $\tau$. 
\hfill $\diamond$
\ere

Let $(G, \piG)$ now be any connected
Poisson Lie group with Lie bialgebra $(\g, \delta_\g)$, and assume that $r \in \g \otimes \g$ is a quasitriangular $r$-matrix for $(\g, \delta_\g)$. Let $Y$ be a manifold with a right $G$-action $\sigma: Y \times G \to Y$,
and assume that the stabilizer subalgebra of $\g$ at every $y \in Y$ is coisotropic with 
respect to the symmetric part of $r$. 
By \leref{le-sigma-pi} and \reref{re-s}, $\sigma(r)$ is a Poisson structure on $Y$, where $\sigma: \g \to \V^1(Y)$ also denotes the right Lie algebra action induced by $\sigma$. 

Recall from $\S$\ref{subsec-quasi-triangle} the pair of dual Lie subalgebras
$((\f_-, \delta_\g|_{\f_-}), (\f_+, -\delta_\g|_{\f_+}))$. 
Let again 
$F_-$ and $F_+$ be the connected subgroups of $G$ with Lie algebras $\f_-$ and $\f_+$ respectively, so  
$(F_-, \,\piG|_{\sF_-})$ and $(F_+, \, -\piG|_{\sF_+})$
form a pair of dual  Poisson Lie groups. 
Restricting the action $\sigma$ of $G$ on $Y$ to actions of $F_\pm$ on $Y$, one then has the action 
groupoids
\[
Y \times F_- \rightrightarrows \,Y \hs \mbox{and} \hs 
Y \times F_+ \rightrightarrows \,Y.
\]
Let $\{x_i\}_{i=1}^n$ be a basis of $\f_-$ and $\{\xi_i\}_{i=1}^n$ the dual basis of $\f_+$ with respect to the
pairing $\lara_{(\f_-, \f_+)}$ between $\f_-$ and $\f_+$ given in \eqref{f_pm-pairing}.

\bco{co-F-pm}
With the notation as above, and let
\begin{align}\lb{eq-pi-p}
\pi_{\sY \times \sF_-} &= (-\sigma(r)) \times_{(\sigma|_{\f_+}, \; \lam_-)} \piG|_{\sF_-} =
(-\sigma(r), \; 0) + (0, \, \piG|_{\sF_-}) -\sum_{i=1}^n (\sigma(\xi_i), \;0) \wedge (0, \; x_i^R),\\
\lb{eq-pi-m}
\pi_{\sY \times \sF_+} &= \sigma(r) \times_{(\sigma|_{\f_-}, \; \lam_+)}(- \piG|_{\sF_+}) =
(\sigma(r), \; 0) + (0, \, -\piG|_{\sF_+}) -\sum_{i=1}^n (\sigma(x_i), \;0) \wedge (0, \; \xi_i^R).
\end{align}
Then 
$(Y\times F_- \rightrightarrows Y, \; \pi_{\sY \times \sF_-})$ and 
$(Y \times F_+ \rightrightarrows Y, \;  \pi_{\sY \times \sF_+})$ form a pair of dual Poisson groupoids.
\eco

\begin{proof} 
Let $\d_{\f_-}$ be the double Lie algebra of $(\f_-, \delta_\g|_{\f_-})$. Then 
$\sigma \circ q: \d_{\f_-} \to \V^1(Y)$ is a Lie algebra homomorphism, where $q: \d_{\f_-} \to \g$ is the
Lie algebra homomorphism given in \eqref{eq-q}. By \reref{re-df-g}, $q(r_{\d_{\f_-}}) = r$. Thus
$(\sigma \circ q)(r_{\d_{\f_-}}) = \sigma(r)$ is a Poisson structure on $Y$. \coref{co-F-pm} now 
follows by applying \coref{co-GG} to the pair of dual Poisson Lie groups 
$(F_-, \,\piG|_{\sF_-})$ and $(F_+, \, -\piG|_{\sF_+})$.
\end{proof}

\bre{re-F-pm}
The Lie algebra action $\sigma \circ q: \d_{\f_-} \to \V^1(Y)$ of $\d_{\f_-}$ on $Y$ gives rise to 
the {\it action Courant algebroid} over $Y$ as defined in \cite{Li-Mein}, with two transversal Dirac structures
defined by the splitting $\d_{\f_-} = \f_- + \f_+$. The pair of dual Poisson groupoids in 
\coref{co-F-pm} then have the two transversal Dirac structures as their Lie bialgebroids.
\hfill $\diamond$
\ere

\sectionnew{Review on standard complex semisimple Poisson Lie groups}\lb{sec-review}

\subsection{The standard complex semisimple Poisson Lie group $(G, \pist)$}\lb{subsec-standard-G}
For the rest of the paper,  let $G$ be a connected complex semisimple Lie group with Lie algebra $\g$. We 
recall the so-called {\it standard multiplicative Poisson structures} on $G$ and refer to
\cite{etingof-schiffmann, Lu-Mou:mixed, Lu-Mou:flags} for details.

Fix a pair $(B, B_-)$ of opposite Borel subgroups of $G$ and a non-degenerate symmetric ad-invariant bilinear form $\lara_\g$ on $\g$, and let $T = B \cap B_-$.  Denote the Lie algebras of $B, B_-$ and $T$ by  $\b, \b_-$ and $\h$ respectively. Let $\g = \h + \sum_{\alpha \in \Delta} \g_\alpha$ be the root decomposition of $\g$ with respect to $\h$, and let $\Delta_+ \subset \h^*$ be the set of positive roots with respect to $\b$. We will also write $\alpha > 0$ for $\alpha \in \Delta_+$.
 Let $\n =  \sum_{\alpha \in \Delta_+} \g_\alpha$, $\n_- =  \sum_{\alpha \in \Delta_+} \g_{-\alpha}$, and let $N$, $N_-$ be the connected subgroups of $G$ with respective Lie algebras $\n$ and $\n_-$. For each $\al >0$, let $E_\alpha \in \g_\alpha$
 and $E_{-\alpha} \in \g_{-\alpha}$ be such that $\la E_\alpha, E_{-\alpha}\ra_\g = 1$. Denote by $\lara$ the bilinear form on both $\h$ and $\h^*$ induced by
$\lara_\g$, and let  $\{h_i\}_{i=1}^{r}$, $r = \dim \h$, be a basis of $\h$ such that $\la h_i, h_j\ra = \delta_{ij}$. The {\it standard quasitriangular $r$-matrix} associated to the choice of the triple $(\b, \b_-, \lara_\g)$ is the element 
\begin{equation}\lb{eq-rst}
\rst = \frac{1}{2} \sum_{i=1}^{r} h_i \otimes h_i + \sum_{\alpha >0} E_{-\al} \otimes E_\al \in \g \otimes \g. 
\end{equation}
The bivector field on $G$ defined by (see notation in $\S$\ref{subsec-nota-intro})
\begin{equation}\lb{eq-pist}
\pist = \rst^L - \rst^R 
\end{equation}
is a multiplicative Poisson structure on $G$, 
and $(G, \pist)$ is called a {\it standard semisimple Poisson Lie group}.
The Lie bialgebra  of $(G, \pist)$ is $(\g, \deltast)$, where $\delta_{{\rm st}}(x) = {\rm ad}_x\rst$ for $x \in \g$. 
In the notation of $\S$\ref{subsec-quasi-triangle}, one has
$$
\Im (r_{{\rm st}}^\sharp) = \b, \hs \text{and} \hs \Im((r_{{\rm st}}^{21})^\sharp) = \b_-.
$$
Thus $B$ and $B_-$ are Poisson Lie subgroups of $(G, \pist)$. Denoting the restrictions of $\pist$ to $B$ and to $B_-$ by the same symbol, the pair $((B_-, \pist), \; (B, -\pist))$ is then a pair of dual  Poisson Lie groups, with the 
pairing $\lara_{(\b_-, \b)}$ in \eqref{f_pm-pairing} given explicitly by 
\begin{equation}\lb{eq-bb-pairing-0}
\la x_- + x_0, \;\, y_+ + y_0\ra_{(\b_-, \b)} = \la x_-,\; y_+\ra_\g + 2\la x_0, \;y_0\ra_\g, \hs x_- \in \n_-,\, 
\, x_0, \,y_0 \in \h, \, y_+ \in \n_.
\end{equation} 
A basis for $\b_-$ and its dual basis for $\b_+$ with respect to the pairing $\lara_{(\b_-, \b)}$ are now given by
\begin{equation}\label{eq-basis}
\{h_i/\sqrt{2}\}_{i=1}^r \cup \{E_{-\alpha}\}_{\alpha > 0} \subset \b_-
\hs \mbox{and} \hs
\{h_i/\sqrt{2}\}_{i=1}^r \cup \{E_{\al}\}_{\alpha > 0}
\subset \b.
\end{equation}

The Poisson structure $\pist$ is invariant under the action of $T$
on $G$ by left or right multiplication. Let $W = N_G(T)/T$ be the Weyl
group of $(G, T)$, where $N_G(T)$ is the normalizer subgroup of $T$ in $G$.
For $u,v \in W$, the {\it double Bruhat cell} (see \cite{Fomin-Zelevinsky:double})
\[
G^{u,v} = BuB \cap B_-vB_-
\]
is non-empty, and $\dim G^{u, v} = l(u) + l(v) + r$, where $l$ is the length function on $W$ and
recall that $r = \dim \h$. 
It is well known \cite{hodges, reshe-4} that the $T$-leaves of $(G, \pist)$ are precisely the double Bruhat cells in $G$. In particular, for each $v \in W$, both $BvB$ and
$B_-vB_-$ are Poisson submanifolds of $G$ with respect to $\pist$.

\subsection{The Drinfeld double and the dressing vector fields of $(G, \pist)$} \lb{subsec-double}
The double Lie algebra $(\d, \lara_\d)$ of the Lie bialgebra $(\g, \delta_{\rm st})$ can be identified with the 
quadratic Lie algebra $(\g \oplus \g, \; \lara_{\gog})$, where $\gog$ has the direct product Lie algebra structure, 
the invariant bilinear form 
$\lara_{\gog}$ is defined by
\[
\la x_1 + y_1,  \, \,  x_2 + y_2 \ra_{\gog} = \la x_1, \,x_2 \ra_\g - \la y_1, \,y_2\ra_\g, 
\hs x_1, \,x_2, \,y_1, \,y_2 \in \g,
\]
and  $\g$ is identified with $\g_{{\scriptscriptstyle \Delta}} = \{(x, x): \, x \in \g\}$ and $\g^*$  with
\begin{equation}\lb{eq-gst}
\gst =  \{(x_+ + x_0, \; -x_0 + x_-): \, 
x_+ \in \n, x_- \in \n^-, x_0 \in \h\}
\end{equation}
(see \cite{chari-pressley, etingof-schiffmann, Lu-Mou:flags}). Let ${r}^{(2)}_{\rm st} \in (\gog) \otimes (\gog)$ be the $r$-matrix on $\gog$ as the double Lie algebra of 
$(\g, \delta_{\rm st})$ (see \exref{exa-double-quasitriangle}), and let
\[
\Pist = \left(r_{\rm st}^{(2)}\right)^L - \left(r_{\rm st}^{(2)}\right)^R
\]
be the corresponding multiplicative Poisson structure on $G \times G$. Then the Poisson Lie group 
$(G \times G, \; \Pist)$ is
a Drinfeld double of $(G, \pist)$, and the diagonal embedding 
\begin{equation}\lb{diag-emb}
(G, \;\pist) \hookrightarrow (G \times G, \;\Pist), \hs g \mapsto (g,\, g),\hs g \in G, 
\end{equation}
realizes $(G, \,\pist)$ as a Poisson subgroup of $(G \times G, \, \Pist)$.

Let $B_-^{{\rm op}}$ be the Lie group which has the same underlying manifold as $B_-$, but with the opposite group structure. 
Then
$$
(\tilde{B}_-, \pi_{\scriptscriptstyle{\tilde{B}_-}}) = (B_- \times B_-^{{\rm op}}, \;\pist \times \pist) \hs \text{and} \hs (\tilde{B}, \,\pi_{\scriptscriptstyle{\tilde{B}}}) = (B \times B, \;(-\pist) \times \pist)
$$
form a pair of dual Poisson Lie groups. Consider the respective right and left Poisson actions 
\begin{align*}
&\rho:\;\; (G, \,\pist) \times (\tilde{B}, \,\pi_{\scriptscriptstyle{\tilde{B}}}) \lrw (G, \,\pist),  \hs \rho(g, \,(b_1, b_2)) = b_1^{-1}gb_2, \hs g \in G, \; b_1, \,b_2 \in B,   \\
&\lam: \;\;(\tilde{B}_-, \, \pi_{\scriptscriptstyle{\tilde{B}_-}}) \times (G, \,\pist) \lrw (G, \,\pist),  \hs 
\lam((b_{-1}, \,b_{-2}), \, g) = b_{-1}gb_{-2}, \hs g \in G, \; b_{-1}, \,b_{-2} \in B_-.
\end{align*}
It is proved in \cite[$\S$6.2 and $\S$8]{Lu-Mou:mixed} that $\Pist$ is a mixed product Poisson structure on $G \times G$.
Namely
\begin{equation} \lb{big-pist}
\Pist = \pist \times_{(\rho, \lam)} \pist. 
\end{equation}

We now present some explicit formulas for the dressing vector fields on $(G, \pist)$
which will be used in the proof of \leref{le-Sigma-uv}.
Let $p_\g: \gog \to \g$ the projection to $\g\cong \gdia$ with respect to the splitting
$\gog = \gdia + \gst$. Note that for any $x \in \g$, writing $x = [x]_- + [x]_0 + [x]_+$
with $[x]_- \in \n^-, [x]_0 \in \h$ and
$[x]_+ \in \n$, one has
\begin{equation}\lb{eq-proj-0}
p_\g(0, x) = \frac{1}{2} [x]_0 +[x]_+  \in \b, \hs p_\g(x, 0) = \frac{1}{2} [x]_0+[x]_-\in \b^-.
\end{equation}
Thus for $\eta \in \n$, the dressing vector field
${\bf d}(\eta, 0)$ at $g \in G$ is given by 
\begin{align}\lb{eq-rho-bb-plus}
{\bf d}(\eta, 0)(g) &=  -l_g p_\g \Ad_{(g^{-1}, g^{-1})} (\eta, 0)  = 
-l_g \left(\frac{1}{2}[\Ad_{g^{-1}}\eta]_0 + [\Ad_{g^{-1}}\eta]_-\right)\\
\nonumber
&=-r_g\eta +l_g\left(\frac{1}{2}[\Ad_{g^{-1}}\eta]_0 + [\Ad_{g^{-1}}\eta]_+\right) \in T_g(gB^-)\cap T_g(BgB),
\end{align}
Similarly, for $\eta \in \n^-$, and $x \in \h$, one has 
\begin{align}\lb{eq-rho-bb-minus}
{\bf d}(0, \eta)(g) & = 
-l_g \left(\frac{1}{2}[\Ad_{g^{-1}}\eta]_0 + [\Ad_{g^{-1}}\eta]_+\right) \\
\nonumber
&=-r_g\eta +l_g\left(\frac{1}{2}[\Ad_{g^{-1}}\eta]_0 + [\Ad_{g^{-1}}\eta]_-\right) \in T_g(gB)\cap T_g(B^-gB^-),\\
\lb{eq-rho-bb-cartan}
{\bf d}(x, -x)(g) &=  l_g \left([\Ad_{g^{-1}}x]_+ - [\Ad_{g^{-1}}x]_-\right)
=r_gx -l_g\left([\Ad_{g^{-1}}x]_0 + 2[\Ad_{g^{-1}}x]_-\right) \\
\nonumber
&=-r_g x +l_g\left([\Ad_{g^{-1}}x]_0 + 2[\Ad_{g^{-1}}x]_+\right) 
\in T_g(TgB^-)\cap T_g(TgB).
\end{align}

\bre{re-also-1}
Note that it also follows from \eqref{eq-rho-bb-plus}, \eqref{eq-rho-bb-minus},  and \eqref{eq-rho-bb-cartan}
that all the $(B, B)$-double cosets and all the $(B_-, B_-)$-double cosets are Poisson submanifold of $(G, \pist)$. 
\hfill $\diamond$
\ere

\subsection{Weak Poisson pairs}\lb{subsec-weak-pairs}
Consider the natural projections
\begin{equation}\lb{eq-varpi-pm}
\varpi: \;\; G \lrw G/B, \;\; g \longmapsto g_\cdot B, \hs \varpi_-:\;\; G \lrw B_-\backslash G, \;\; g \longmapsto {B_-}_\cdot
g, \hs g \in G.
\end{equation}
As both $B$ and $B_-$ are Poisson Lie subgroups of $(G, \pist)$, 
\begin{equation}\lb{eq-pi-12}
\pi_1 \stackrel{{\rm def}}{=}\varpi(\pist) \hs \mbox{and} \hs \pi_{-1} \stackrel{{\rm def}}{=} \varpi_-(\pist)
\end{equation}
are now well-defined Poisson structures on $G/B$ and on $B_-\backslash G$, respectively. 
The Poisson structure $\pi_1$ is invariant under the action of $T$ on $G/B$ by left multiplication, and it is proven in \cite{Go-Ya} that the $T$-leaves of $\pi_1$ are precisely the so-called {\it open Richardson varieties}, i.e non-empty intersections $(BuB/B) \cap (B_-wB/B)$, where $u,w \in W$. In particular, every {\it Bruhat cell} $BuB/B$, for $u \in W$, is a Poisson subvariety of $(G/B, \,\pi_1)$. Similarly, every Bruhat cell $B_-\backslash B_-uB_-$, for $u \in W$,
is a Poisson subvariety of $(B_-\backslash G, \, \pi_{-1})$.

\bde{de-weak-pair} \cite[$\S$8.6]{Lu-Mou:mixed} Two Poisson maps
$\rho_\sY: (X, \piX) \to (Y, \piY)$ and $\rho_{\sZ}: (X, \piX) \to (Z, \piZ)$ are said to form a
{\it weak Poisson pair} if the map
\[
(\rho_\sY, \; \rho_\sZ): \;\;\; (X, \, \piX) \lrw (Y \times Z, \; \piY \times \piZ),\;\;
(y, \, z) \longmapsto (\rho_\sY(y), \; \rho_\sZ(z)), \hs y \in Y, \, z \in Z,
\]
is a Poisson map.
\ede

The following 
\leref{le-weak-pair} is a special case of a fact proved in \cite[$\S$8.6]{Lu-Mou:mixed}, but
for the convenience of the reader, we give a
proof which is much simpler in our special case. 

\ble{le-weak-pair}
The two Poisson maps 
\[
\varpi: \;\; (G, \,\pist) \lrw (G/B, \;\pi_1) \hs \mbox{and} \hs
\varpi_-: \;\; (G, \,\pist) \lrw (B_-\backslash G, \;\pi_{-1}) 
\]
form a weak Poisson pair. Consequently, for $u, v \in W$ and for any symplectic leaf $\Sigma^{u, v} \subset \Guv$,
one has the weak Poisson pairs
\begin{align*}
&\varpi|_{\sG^{u, v}}: \;\; (\Guv, \,\pist) \lrw (BuB/B, \;\pi_1) \hs \mbox{and} \hs
\varpi_-|_{\sG^{u, v}}: \;\; (\Guv, \,\pist) \lrw (B_-\backslash B_-vB_-, \;\pi_{-1}),\\
&\varpi|_{\sSigma^{u, v}}: \;\; (\Sigma^{u, v}, \,\pist) \lrw (BuB/B, \;\pi_1) \hs \mbox{and} \hs
\varpi_-|_{\sSigma^{u, v}}: \;\; (\Sigma^{u, v}, \,\pist) \lrw (B_-\backslash B_-vB_-, \;\pi_{-1}).
\end{align*}
\ele

\begin{proof}
Consider the projection $\Phi: G \times G \to (G/B) \times ({B_-}\backslash G)$ defined by
\[
\Phi(g_1, \, g_2) = ({g_1}_\cdot B, \; {B_-}_\cdot g_2), \hs g_1, \, g_2 \in G.
\]
Using \eqref{big-pist} to write $\Pist = (\pist, \, 0) + (0, \, \pist) + \pi_{\rm mix}$, it follows from 
the definition of the mixed term $\pi_{\rm mix}$ that $\Phi(\pi_{\rm mix}) = 0$. Thus
\[
\Phi: \;\; (G \times G, \; \Pist) \lrw ((G/B) \times (B_-\backslash G), \; \pi_1 \times \pi_{-1})
\]
is Poisson. As the diagonal embedding $(G, \pist) \hookrightarrow (G \times G, \, \Pist)$ is Poisson,
$\varpi$ and $\varpi_-$ form a weak Poisson pair. As $\Guv$ or any symplectic leaf in $\Guv$ are Poisson submanifolds of $(G, \pist)$, the rest of \leref{le-weak-pair} follows.
\end{proof}

Note that in \deref{de-weak-pair}  we do not require the two maps $\rho_\sY$ and $\rho_\sZ$ in a weak Poisson pair to be surjective nor submersions. The next \leref{le-Sigma-uv} and \leref{le-Guv-uv} say that the Poisson maps in the weak Poisson pairs in \leref{le-weak-pair}, although not necessarily surjective, are all submersions.

\ble{le-Sigma-uv}
For any $u, v\in W$ and any symplectic leaf $\Sigma^{u, v}$ of $\pist$ in $\Guv$,
the maps 
\[
\varpi|_{\sSigma^{u, v}}: \;\;\Sigma^{u, v} \lrw BuB/B \hs \mbox{and} \hs \varpi_-|_{\sSigma^{u, v}}: \;\;
\Sigma^{u, v} \lrw
B_-\backslash B_-vB_-
\]
are submersions. 
\ele

\begin{proof} Let $g \in \Sigma^{u, v}$. By definition,
the value at $g$ of every dressing vector field on $(G, \pist)$ is tangent to $\Sigma^{u, v}$. 
By \eqref{eq-rho-bb-plus} and \eqref{eq-rho-bb-cartan}, the differential
of $\varpi|_{\sSigma^{u, v}}$ at $g$ is a surjective linear
map from $T_g\Sigma^{u, v}$ to $T_{g_\cdot B}(BuB/B)$. Thus 
$\varpi|_{\sSigma^{u, v}}: \Sigma^{u, v} \to BuB/B$ is a submersion. Similarly, 
$\varpi_-|_{\sSigma^{u, v}}: \Sigma^{u, v} \to B_-\backslash B_-vB_-$ is a submersion.
\end{proof}

\bre{re-another} \leref{le-Sigma-uv} implies that for any $u, v \in W$, 
the maps 
\[
\varpi|_{\sG^{u, v}}: \;\; (\Guv, \,\pist) \lrw (BuB/B, \;\pi_1) \hs \mbox{and} \hs
\varpi_-|_{\sG^{u, v}}: \;\; (\Guv, \,\pist) \lrw (B_-\backslash B_-vB_-, \;\pi_{-1})
\]
are also submersions, a fact one can in fact see directly without computing the dressing vector fields.
 Indeed, For any 
$g \in G$ and $x \in \b$, the element 
\[
z_{g, x} \, \stackrel{\rm def}{=} \, 
r_g x -l_g\left(\frac{1}{2}\left(\left[\Ad_{g^{-1}}x\right]_0\right) + \left[\Ad_{g^{-1}}x\right]_+\right)
=l_g \left(\frac{1}{2}\left(\left[\Ad_{g^{-1}}x\right]_0\right) + \left[\Ad_{g^{-1}}x\right]_-\right)
\]
lies in $T_g(BgB \cap B_-gB_-)$ and $\varpi(z_{g, x}) = \varpi(r_g x)$.
It follows that the differential of $\varpi$ restricts to a surjective linear map
from $T_g(BgB \cap B_-gB_-)$ to $T_{g_\cdot B}(Bg_\cdot B)$  for
every $g \in G$. This shows in particular that for  
any $u, v\in W$, the map $\varpi|_{\sG^{u, v}}: G^{u, v} \to BuB/B$ is a submersion. Similarly, one sees
that $\varpi_-|_{\sG^{u, v}}$ is a  submersion.
\hfill $\diamond$
\ere

\ble{le-Guv-uv} 
For any $u, v \in W$ and for any symplectic leaf $\Sigma^{u, v}$ of $\pist$ in $\Guv$, one has 
\begin{align*}
&\varpi(\Sigma^{u, v}) = \varpi(G^{u, v}) = \bigcup_{w \leq u, \, w \leq v} (BuB/B) \cap (B_-wB/B)
\subset BuB/B,\\
&\varpi_-(\Sigma^{u, v}) = \varpi_-(G^{u, v}) = \bigcup_{w \leq u, \, w \leq v}
(B_-\backslash B_-wB) \cap (B_-\backslash B_-vB_-) \subset B_-\backslash B_- v B_-,
\end{align*}
where $\leq$ is the Bruhat order on $W$ defined by the
the choice of $B$.
\ele

\begin{proof}
For $w \in W$, $B_-wB \subset B_-vB_-B$ if and only if $B_-wB \cap B_-vB_- \neq \emptyset$,
which, by \cite[Corollary 1.2]{De}, is equivalent to $w \leq v$.  Thus 
$B_- v B_-B = \bigcup_{w \leq v} B_-wB$. It follows that 
\[
\varpi(\Guv) =  \varpi(BuB) \cap \varpi(B_-vB_-B) = 
\bigcup_{w \leq u, \, w \leq v} (BuB/B) \cap (B_-wB/B).
\]
Since $\Guv = \Sigma^{u, v}T$, one has $\varpi(\Sigma^{u, v}) = \varpi(\Guv)$. The claims on 
$\varpi_-(\Sigma^{u, v})$ and $\varpi_-(\Guv)$ are proved similarly.
\end{proof}

\bre{re-symplectic-pair} For $u, v \in W$ and a symplectic leave $\Sigma^{u, v}$ of $\pist$ in $\Guv$,
the weak Poisson pair  
\[
\varpi|_{\sSigma^{u, v}}: \;\; (\Sigma^{u, v}, \,\pist) \lrw (BuB/B, \;\pi_1) \hs \mbox{and} \hs
\varpi_-|_{\sSigma^{u, v}}: \;\; (\Sigma^{u, v}, \,\pist) \lrw (B_-\backslash B_-vB_-, \;\pi_{-1})
\]
in \leref{le-weak-pair} is in general not a {\it symplectic dual pair} \cite{Wein:local}
which requires the two Poisson maps be surjective submersions and their fibers  mutual symplectic othogonals
of each other.
\hfill $\diamond$
\ere


\sectionnew{The double Bruhat cells $\Gvv$ as Poisson groupoids} \lb{sec-ss-PL}
Let the notation be as in $\S$\ref{sec-review}. 
In this section, we apply the results in $\S$\ref{sec-gpoids-algoids} to the  Poisson
Lie group $(G, \pist)$ to construct an action Poisson groupoid $((G/B) \times B_-, \pi)$ over $(G/B, \pi_1)$. For 
$v \in W$, the choice of a representative $\bv$ of $v$ in $N_\sG(T)$ is used to identify $(\Gvv, \pist)$ 
with a Poisson subgroupoid of $((G/B) \times B_-, \pi)$ through a Poisson embedding
$I_\bv: (B_-vB_-, \pist) \hookrightarrow ((G/B) \times B_-, \pi)$. 

\subsection{The action Poisson groupoid $((G/B) \times B_-, \pi)$ over $(G/B, \pi_1)$} \lb{subsec-BGB}
Let $G$ act on the flag variety $G/B$ from the {\it right} by 
\[
(G/B) \times G \lrw G/B, \;\; (g_\cdot B, \; g_1) \longmapsto g_1^{-1}g_\cdot B,\hs g, \, g_1 \in G,
\]
and let $\sigma: \g \to \V^1(G/B)$ be the induced {\it right} Lie algebra action of $\g$ on $G/B$ given by
\begin{equation}\lb{eq-sigma-GB}
\sigma(x) = -\varpi(x^R), \hs \mbox{or}\hs \sigma(x)(g_\cdot B) = \frac{d}{dt}|_{t=0} \exp(-tx)g_\cdot B
\hs x \in \g, \, g \in G,
\end{equation}
where recall that $\varpi: G \to G/B$ is the projection.
Restricting the $G$-action on $G/B$ to one of $B_-$ on $G/B$, one then has the action groupoid 
$(G/B) \times B_- \rightrightarrows G/B$, with the source map $\theta$, the target map $\tau$,
the groupoid multiplication $\mu$,
the inverse map $\iota$, and the identity bisection $\epsilon$ respectively given by
\begin{align}\lb{eq-maps-0}
&\theta(g.B, \; b_-) = g.B,  \hs \;\;\;\tau(g.B, \; b_-) = (b_-^{-1}g).B,  \\
\lb{eq-maps-15}
&\mu(g_\cdot B, \; b_-, \; b_-^{-1}g_\cdot B, \; b_-^\prime) = (g_\cdot B, \; b_-b_-^\prime),\\
\lb{eq-maps-1}
&\iota(g.B, \; b_-) = (b_-^{-1}g.B, \; b_-^{-1}),\hs \;\;\;\epsilon(g_\cdot B) = (g_\cdot B, \, e),\hs
b_-, \, b_-^\prime \in B_-, \; g \in G.
\end{align}
Consider the Poisson structure $\pi_1 =\varpi(\pist)$ on $G/B$. As $\pist = r_{\rm st}^L - r_{\rm st}^R$ and 
$\varpi(r_{\rm st}^L) = 0$, one has
\begin{equation}\lb{formula_pi_1}
\pi_1 =- \varpi(r_{\rm st}^R) = -\sigma(\rst). 
\end{equation} 
Let $\lam_-: \b_- \to \V^1(B_-)$ be given by $\lam_-(x) = x^R$ for $x \in \b_-$.
As $\sigma|_\b: \b \to \V^1(G/B)$ is a right Poisson action of the Lie bialgebra $(\b, -\delta_{\rm st}|_\b)$ on
$(G/B, \pi_1)$, one has the mixed product Poisson structure $\pi$ on $(G/B) \times B_-$ given by
\begin{equation} \lb{defn-pi-gpoid}
\pi = \pi_1  \times_{(\sigma|_\b, \,\lam_-)} \pist = (\pi_1, \, 0) + (0, \, \pist)
-\sum_{i=1}^n (\sigma(\xi_i), \;0) \wedge (0, \; x_i^R),
\end{equation}
where $\{x_i\}_{i=1}^n$ is any basis of $\b_-$ and $\{\xi_i\}_{i=1}^n$ the dual basis of $\b$ 
with respect to the with the 
pairing $\lara_{(\b_-, \b)}$ between $\b_-$ and $\b$ given in \eqref{f_pm-pairing}. 
By \coref{co-F-pm} and \eqref{formula_pi_1}, 
\[
((G/B) \times B_- \rightrightarrows G/B, \;\; \pi)
\] 
is an action Poisson groupoid over the Poisson manifold $(G/B, \pi_1)$. Note that the bases for $\b_-$ and $\b$ in \eqref{defn-pi-gpoid} can be taken to be the ones in  
\eqref{eq-bb-pairing-0}. 

\subsection{The Poisson embedding of $(B_- v B_-, \pist)$ into $((G/B) \times B_-, \, \pi)$}\lb{subsec-embedding}
Recall that $N_\sG(T)$ is the normalizer subgroup of $T$ in $G$. In this section, fix 
$v \in W$ and let $\bv \in N_\sG(T)$ be any representative of $v$ in $N_\sG(T)$. Let 
\begin{equation} \lb{defn-E_u}
C_\bv = N\bv \cap \bv N_- \subset G.
\end{equation}
It is well known that the multiplication maps 
\begin{align*}
&C_\bv \times B\, \lrw \, BvB, \; \;\;(c, \, b) \longmapsto cb, \hs c \in C_\bv, \, b \in B,\\
&B_- \times C_\bv \, \lrw \, B_-vB_-, \;\;\; (b_-,\, c) \longmapsto b_-c, \hs b_- \in B_-, \; c \in C_{\bv},
\end{align*}
are algebraic isomorphisms. 
Consider now the embedding
\begin{equation}\lb{eq-Iv-0}
I_\bv: \;\; B_-vB_- \lrw  (G/B) \times B_-, \;\; I_\bv(b_-c) = (b_-c_\cdot B, \; b_-), 
\hs b_- \in B_-, \, c \in C_\bv.
\end{equation}
The goal of $\S$\ref{subsec-embedding} is to prove the following \prref{pr-embed}.

\bpr{pr-embed}
The embedding $I_\bv: (B_-vB_-, \,\pist) \to ((G/B) \times B_-, \, \pi)$ is Poisson. 
\epr

To prepare for the proof of \prref{pr-embed}, we first prove some properties of $C_\bv$. 

\ble{le-Cv-0}
The submanifold $C_\bv$ of $G$ is coisotropic with respect to the Poisson structure $\pist$. 
\ele

\begin{proof} Consider first the 
subgroup $N_v =N \cap (\bv N_- \bv^{-1})$ with Lie algebra $\n_v = \n \cap \Ad_{\bv} \n_-$. We first show that
$N_v \subset G$ is coisotropic with respect to $\pist$. 
With $\g^* \cong \gst$, where the pairing between $\g \cong \gdia$ and
$\gst$ is via the bilinear form $\lara_{\gog}$ on $\gog$, the annihilator subspace 
$\n_v^0 = \{\xi \in \g^*: \xi|_{\n_v} = 0\}$ of $\n_v$ in $\gst$ is 
\[
\{(x_+ + x_0, \; -x_0 + x_-): \; x_+ \in \n, \, x_0 \in \h, \, x_- \in \n_- \cap \Ad_{\bv} \n\},
\]
which is a Lie subalgebra of $\gst$. It follows \cite{lu-we:Poi, STS2} that $N_v$
is a coisotropic subgroup of $(G, \pist)$. 

Let $c \in C_\bv$ and write $c = n\bv$, where $n \in N_v$. By the multiplicativity of $\pist$, one has
\[
\pist(c) = \pist(n\bv) = l_{n} \pist(\bv) + r_{\bv} \pist(n).
\]
As $N_v$ is coisotropic with respect to $\pist$, $\pist(n) \in (T_{n}G) \wedge (T_{n} N_v)$, so
$r_{\bv} \pist(n) \in (T_c G) \wedge (T_c C_\bv)$.   On the other hand, it is
easy to see that
\begin{equation}\lb{eq-pist-v}
\pist(\bv) = -r_\bv\left(\sum_{\alpha > 0, v^{-1}\alpha < 0} 
E_{-\alpha} \wedge E_\alpha\right).
\end{equation}
It follows that $l_n \pist(\bv) \in (T_cG) \wedge (T_cC_\bv)$. Thus 
$C_\bv$ is a coisotropic submanifold of $(G, \pist)$.
\end{proof}

\ble{proj-to-B-poiss-0}
The map 
$$
q_\bv: \;\;(B_-vB_-, \, \pist) \to (B_-, \, \pist), \hs q_\bv(b_-c) = b_-, \hs b_- \in B_-, \; c \in C_\bv,
$$
is Poisson. 
\ele

\begin{proof}
(See also \cite[Theorem 3.1]{GSV}) Let $b_- \in B_-$ and $c \in C_\bv$. By the multiplicativity of $\pist$, one has $\pist(b_-c) = l_{b_-}\pist(c) + r_c \pist(b_-)$. As $C_\bv$ is a coisotropic submanifold of $(B_-vB_-, \pist)$, one has $\pist(c) \in T_cC_\bv \wedge T_c(B_-vB_-)$. As $q_\bv(l_{b_-}T_cC_\bv) = 0$, one has $q_\bv l_{b_-}\pist(c)=0$. Using the fact that $\pist(b_-) \in \wedge^2 T_{b_-} B_-$, one sees that 
\[
q_\bv(\pist(b_-c)) = (q_\bv r_c) (\pist(b_-)) = \pist(b_-).
\] 
\end{proof}

\noindent
{\it Proof of \prref{pr-embed}}: Let
$(B_-vB_-)_{\rm diag} = \{(g, \,g): \,g \in B_-vB_-\}$. 
Then $I_\bv$ is the restriction to $(B_-vB_-)_{\rm diag}  \subset G \times (B_-vB_-)$ of the map
\[
K_\bv: \;\;G \times (B_-vB_-) \lrw (G/B) \times B_-, \;\; (g, \, b_-c) \longmapsto (g_\cdot B, \; b_-), 
\hs g \in G, \, b_- \in B_-, \, c \in C_\bv.
\]
By $\S$\ref{subsec-double} and in particular \eqref{big-pist}, both $(B_-vB_-)_{\rm diag}$ and
$G \times (B_-vB_-)$ are Poisson submanifolds of
$G \times G$ with respect to the Poisson structure $\Pist$. It is thus enough to show that 
\[
K_\bv: \;\;(G \times (B_-vB_-), \; \Pist) \lrw ((G/B) \times B_-, \; \pi)
\]
is Poisson.
Let again $(x_i)_{j=1}^n$ be any basis of $\b_-$ and  $(\xi_i)_{j=1}^n$ the basis of $\b$ dual to $(x_i)_{i=1}^n$ under the 
pairing $\lara_{(\b_-, \b)}$ in \eqref{eq-bb-pairing-0}.
By \eqref{big-pist}, one has $\Pist = (\pist, 0) + (0, \pist) + \mu_1 + \mu_2$, where 
$$
\mu_1 = \sum_{i=1}^n (\xi_i^R, 0) \wedge (0, x_i^R), \hs \text{and} \hs \mu_2 = -\sum_{i=1}^n (\xi_i^L, 0) \wedge (0, x_i^L).
$$
By the definition of $\pi_1$, $K_\bv(\pist, 0) = (\pi_1, 0)$. By \leref{proj-to-B-poiss-0}, $K_\bv(0, \pist) = (0, \pist)$.
Since for any $\xi \in \b$, the vector field $\xi^L$ on $G$ vanishes when projected to $G/B$, one has $K_\bv(\mu_2) = 0$. 
It is also clear from the definitions that $K_\bv (\mu_1)$ coincides with the mixed term of $\pi$. 
Thus $K_\bv$ is Poisson.

This finishes the proof of \prref{pr-embed}.

\bre{BvB-mixed_prod} {\bf (The Poisson structure $\pist$ on $B_-vB_-$ as a mixed product)} Define
\[
\Psi: \;\;\; (G/B) \times B_- \lrw B_- \times (G/B), \;\;\;\Psi(g_\cdot B, \; b_-) = (b_-^{-1}, \; g_\cdot B),
\]
and consider the Poisson structure $\pi^\prime = -\Psi(\pi)$ on $B_- \times (G/B)$. It is easy to see that 
\[
\pi' = \pist \times_{(\rho_-, \, \lambda_+)}  (-\pi_1),
\]
where $\rho_-$ and $\lam_+$ denote the Poisson Lie group actions as well as the induced Lie bialgebra actions, respectively given 
by
\begin{align*}
&(B_-, \, \pist) \times (B_-, \, \pist) \lrw (B_-, \pist), \;\;\; (b_-, \, b_-^\prime) \longmapsto
b_- b_-^\prime, \hs b_-, \, b_-^\prime \in B_-,\\
&(B_+, -\pist) \times (G/B, -\pi_1) \lrw (G/B, -\pi_1), \;\; (b, \; g_\cdot B) \longmapsto bg_\cdot B, \; \; \; b\in B,
\, g \in G.
\end{align*}
One then has the Poisson embedding
\begin{equation}\lb{eq-BvB-m}
\Psi \circ\iota \circ I_\bv: \;\; (B_-vB_-, \,\pist) \lrw (B_- \times (G/B), \, \pi^\prime),
\;\; b_-c \longmapsto  (b_-, \, c_\cdot B),
\;\;\;\; b_- \in B_-, \, c \in C_\bv,
\end{equation}
where $\iota$ is the inverse map of the Poisson groupoid
$((G/B) \times B_-\rightrightarrows G/B, \,\pi))$.
Note the image of $B_-vB_-$ under $\Psi \circ\iota \circ I_\bv$ is the Poisson submanifold
$B_- \times (BvB)/B$ of $(B_- \times (G/B), \;\pi^\prime)$. We have thus 
identified the restriction of $\pist$ to $B_-v B_-$ as the mixed product Poisson structure $\pi'$ 
on the product manifold $B_- \times (BvB/B)$ via the map in \eqref{eq-BvB-m}.  \hfill $\diamond$
\ere

\bre{re-BvB-match} 
Consider also the Poisson embedding
\[
J_\bv \, \stackrel{{\rm def}}{=}\, \iota \circ I_\bv: \; (B_-vB_-, -\pist) \lrw
((G/B) \times B_-, \, \pi),\; J_\bv(b_-c) = (c_\cdot B, \, b_-^{-1}), \;\;\; b_- \in B_-, \, c \in C_\bv.
\]
Then $J_\bv(B_-vB_-) = (BvB/B) \times B_-$.  As $v$ runs over $W$, one has the respective disjoint unions
\[
G = \bigsqcup_{v \in W} B_-vB_- \hs \mbox{and} \hs (G/B) \times B_- = \bigsqcup_{v \in W} (BvB/B) \times B_-
\]
of the Poisson varieties $(G, -\pist)$ and $((G/B) \times B_-, \pi)$ into Poisson subvarieties, together with
piecewise Poisson isomorphisms $\{J_\bv: v \in W\}$, but these  piecewise Poisson isomorphisms do not patch together to define a smooth map from $G$ to $(G/B) \times B_-$.
\hfill $\diamond$
\ere

\bex{broken-isom-exa} 
Let $G = SL(2, \Cset)$ and let $B$ and $B_-$ be the subgroups of $G$
consisting of upper and lower triangular matrices
respectively. Let $s \in W$ be the non-trivial element, so that 
$$
B_-sB_- = \left\{ \left(\begin{array}{cc}a & b \\c & d\end{array}\right) : ad - bc = 1, \; b \neq 0 \right\}. 
$$
Identify the flag variety $G/B$ with the complex projective space $\cp$ via $\left(\begin{array}{cc}a & b \\c & d\end{array}\right).B \mapsto [a,c]$. 
For $\bar{s}= \!\left(\!\begin{array}{cc}0 & -1 \\1 & 0\end{array}\!\!\right)$,
the map $J_{\bar{s}}: B_-sB_- \to \cp \times B_-$ is given by 
\[
J_{\bar{s}} \left(\begin{array}{cc}a & b \\c & d\end{array}\right) = \left([a, \, -b], \;  \left(\begin{array}{cc} -b^{-1} & 0 \\d & -b \end{array}\right)\right),
\]
which does not extend to a smooth map from $G$ to $\cp \times B_-$.
\hfill $\diamond$
\eex

\subsection{Poisson embeddings of $(\Guv, \,\pist)$ into $((G/B) \times B_-, \, \pi)$}\lb{subsec-Guv-embedding}
Recall that $\theta$ and $\tau$
 are respectively the source and target maps of
the action groupoid $(G/B) \times B_-$ over $G/B$, and note that 
the image of $B_-vB_-$ under the embedding $I_\bv$ is 
\[
I_\bv(B_-vB_-) = \tau^{-1}(BvB/B) = \iota((BvB/B) \times B_-).
\]
For $u \in W$, restricting $I_\bv$ to $\Guv  =BuB \cap B_-vB_- \subset B_-vB_-$, one has the embedding
\begin{equation}\lb{eq-I-guv-GB}
I_\bv|_{\sG^{u, v}}: \;\;\; \Guv \hookrightarrow (G/B) \times B_-.
\end{equation}
For $u, v \in W$, set 
\begin{equation}\lb{eq-Fuv}
F^{u, v} \; \stackrel{\rm def}{=} \; \theta^{-1}(BuB/B) \cap \tau^{-1}(BvB/B) \subset (G/B) \times B_-, \hs u,v \in W. 
\end{equation}
It is clear
from the definitions that
\begin{equation}\lb{eq-Jv-Guv-0}
I_\bv (\Guv) = F^{u, v}, \hs \hs u \in W.
\end{equation}

Let $T$ act on $(G/B) \times B_-$  via
\begin{equation}\lb{eq-T-actions}
t \cdot (g_\cdot B, \; b_-) = (tg_\cdot B, \; tb_-), \hs t \in T, \, g \in G, \, b_- \in B_-.
\end{equation}

\bpr{T-leaves-pi-0} The mixed product Poisson structure  $\pi$ on $(G/B) \times B_-$ is
invariant under the $T$-action, and its $T$-leaves  
are precisely the intersections $F^{u, v}$, where $u, v \in W$.
\epr

\begin{proof} For each $v \in W$, choose a representative $\bv$ of $v$ in $N_\sG(T)$.
Let $T$ act on $B_-vB_-$ by left translation. Clearly, $I_\bv: B_- v B_- \to (G/B) \times B_-$ is $T$-equivariant.
The statement of \prref{T-leaves-pi-0} now follows from the $T$-equivariant Poisson isomorphisms $I_\bv$, $v \in W$, and the fact
that the $T$-leaves of $\pist$ in $B_-vB_-$ are the $\Guv$'s for $u \in W$. 
\end{proof}

\bre{re-twist} ({\bf The Fomin-Zelevinsky twist map})
Let $u, v \in W$ and let $\bu$ and $\bv$ be any representatives of $u$ and $v$ in $N_\sG(T)$ respectively. 
Recall that the inverse map $\iota$ of the Poisson groupoid $((G/B) \times B_-, \, \pi)$ satisfies
$\iota(\pi) = -\pi$. As $\iota(F^{u, v}) = F^{v, u}$, by \prref{pr-embed},
\begin{equation}\lb{eq-iota-uv}
\iota^{\bu, \bv} \, \stackrel{{\rm def}}{=}\, \left(I_\bu|_{\sG^{v, u}}\right)^{-1} \circ \iota \circ \left(I_\bv|_{\sG^{u, v}}\right): \;\;\; (\Guv, \; \pist) \lrw 
(G^{v, u}, \; \pist)
\end{equation}
is anti-Poisson. Explicitly, the map $\iota^{\bu, \bv}: G^{u, v} \to G^{v, u}$ is given by
\[
\iota^{\bu, \bv}(g) = b_-^{-1} c = c'b^{-1}, \hs \mbox{if} \;\; g = cb= b_-c^\prime \in G^{u, v},\; \,\mbox{where}\;\,
c \in C_\bu, \, b \in B, \, b_- \in B_-, \, c' \in C_\bv,
\]
or, if for $h \in N_-TN$, we write $h = [h]_-[h]_0[h]_+$ with $[h]_- \in N_-, [h]_0 \in T, [h]_+ \in N$, then
\begin{equation}\lb{eq-tau-uv}
\iota^{\bu, \bv}(g) = \left([\bu^{-1}g]_-^{-1} \bu^{-1}g\bv^{-1}[g\bv^{-1}]_+^{-1}\right)^{-1}, \hs g \in \Guv.
\end{equation}
In \cite[$\S$1.5]{Fomin-Zelevinsky:double}, Fomin and Zelevinsky introduced a twist map $G^{u, v} \to G^{u^{-1}, v^{-1}}$
(for certain special ways of choosing $\bu$ and $\bv$). By \eqref{eq-tau-uv},
the Fomin-Zelevinsky twist map
is the composition of $\iota^{\bu, \bv}$ with the group inverse $G \to G, g \mapsto g^{-1}$, of $G$ and with
an involutive automorphism $x \to x^\theta$ of $G$ (see \cite[Formula (1.11)]{Fomin-Zelevinsky:double}, while 
the latter two involutions are easily seen to
be both anti-Poisson with respect to $\pist$. It follows that the 
Fomin-Zelevinsky twist $(G^{u, v}, \pist) \to (G^{u^{-1}, v^{-1}}, \pist)$ is anti-Poisson, a fact already proved
in \cite[Theorem 3.1]{GSV}.
\hfill $\diamond$
\ere

\bre{disjoint-union-dlbe-B-cell-0} Consider the two disjoint union decompositions 
\begin{equation} \lb{disjoint-union-formula-0}
G = \bigsqcup_{u,v \in W} G^{u,v}, \hs \hs (G/B) \times B_- = \bigsqcup_{u,v \in W} F^{u,v}. 
\end{equation}
Let $T$ act on $G$ by left multiplication and on $(G/B) \times B_-$ by  \eqref{eq-T-actions}.
Then the two decompositions in \eqref{disjoint-union-formula-0} are respectively that of $T$-leaves of
$(G, \pist)$ and $((G/B) \times B_-, \pi)$. Any choice $\{\bv \in N_\sG(T): v \in W\}$ gives rise to
piecewise $T$-equivariant Poisson isomorphisms
\[
I_\bv: \;\; (B_- vB_- = \bigsqcup_{u \in W} \Guv, \; \, \pist) \lrw 
(\tau^{-1}(BvB/B) =\bigsqcup_{u \in W} F^{u, v},\; \pi)
\]
but the maps $\{I_\bv: v \in W\}$ do not patch together to define a smooth map from $G$ to $(G/B) \times B_-$.
See also \reref{re-BvB-match}.
\hfill $\diamond$
\ere

\subsection{The double Bruhat cell $\Gvv$ as Poisson groupoids}\lb{subsec-Gvv-groupoid}
Observe that for any $v \in W$, 
\[
F^{v, v} = \theta^{-1}(BvB/B) \cap \tau^{-1}(BvB/B) \subset (G/B) \times B_-
\]
is the subgroupoid of $(G/B) \times B_-\rightrightarrows G/B$ over the subset $BvB/B$ of $G/B$.

\bde{defn-G^bvbv}
For $v \in W$ and any representative $\bv$ of $v$ in $N_\sG(T)$, denote by $G^{\bv, \bv} \rightrightarrows BvB/B$
the double Bruhat cell $G^{v,v}$, equipped with the groupoid structure induced by the isomorphism
$I_\bv: \Gvv \to F^{v, v}$. In details, the groupoid structure is defined as follows: 
for $g = c b = b_-c' \in G^{v,v}$, where $b \in B$, $b_- \in B_-$, and $c, c' \in C_\bv$, 
\begin{align*}
&\text{source map}: \; \theta_\bv(g) = g.B = c.B,   \\
&\text{target map}: \; \tau_\bv(g) = c'.B,    \\
&\text{inverse map}: \; \iota_\bv(g) = c'b^{-1} = b_-^{-1}c,   \\
& \text{identity bisection}: \; \epsilon_\bv(c_\cdot B) = c \in C_\bv \subset G^{v,v}.
\end{align*}
If $h \in G^{v,v}$ is such that $\tau_\bv(g) = \theta_\bv(h)$, so $h =  c'b' = b_-'c''$, with $b' \in B$, $b'_- \in B_-$, and $c'' \in C_\bv$, the groupoid product of $g$ and $h$ is  given by 
\begin{equation} \lb{mult-in-Gvv}
\mu_\bv(g, h) = c bb' = b_-b_-'c''. 
\end{equation}
\ede

The following \thmref{poiss-gpoid-Gvv}, which follows directly from \prref{pr-embed}, is the first main result of this paper. 

\bthm{poiss-gpoid-Gvv}
For any $v \in W$ and $\bv \in N_G(T)$, the pair $(G^{\bv, \bv}, \pist)$ is a Poisson groupoid over the Poisson manifold
$(BvB/B, \pi_1)$. 
\ethm

\begin{proof} 
It is clear that all the structure maps of the groupoid $G^{\bv, \bv} \rightrightarrows BvB/B$ are smooth. 
As $C_\bv \subset G^{\bv, \bv}$, the source map $\theta_\bv$ is surjective. By \leref{le-Sigma-uv}, 
$\theta_\bv$ is a submersion. Thus  $G^{\bv, \bv}$ is a Lie groupoid over
$BvB/B$. As $I_\bv(\Gvv)$ is a Poisson submanifold of $(G/B) \times B_-$ with respect to $\pi$,
$(G^{\bv, \bv}, \pist)$ is a Poisson groupoid over
$(BvB/B, \pi_1)$.
\end{proof}

\bre{all-isomo}
If $\bv, \tilde{v}$ are two representatives of $v \in W$ and if $t \in T$ is such that $\bv = t\tilde{v}$,
then the left translation $l_t: (G^{\tilde{v}, \tilde{v}}, \pist) \to (G^{\bv, \bv}, \pist)$ is a Poisson
groupoid isomorphism covering the Poisson 
isomorphism $l_t: (BvB/B, \pi_1) \to (BvB/B, \pi_1)$. Hence the isomorphism class of $(G^{\bv, \bv}, \pist)$ as a Poisson groupoid is independent of the choice of the representative $\bv$.  \hfill $\diamond$ 
\ere

Recall that $\varpi_-: G \to B_-\backslash G$ is the 
projection, and for each $v \in W$, $B_-\backslash B_-vB_-$ is a Poisson submanifold of $B_-\backslash G$ with respect to the Poisson structure $\pi_{-1} = \varpi_-(\pist)$. For $v \in W$ and any representative $\bv$ of $v$ in $N_\sG(T)$, define
\begin{equation}\lb{eq-Phi-bv}
\Phi_\bv: \;\; B_-\backslash B_-vB_- \lrw BvB/B, \;\;\; {B_-}_\cdot c \longmapsto c_\cdot B, \hs c \in C_\bv.
\end{equation}

\ble{le-Phi-bv} 
For $v \in W$ and any representative $\bv$ of $v$ in $N_\sG(T)$,
\[
\Phi_\bv:  \;\;\;(B_-\backslash B_-vB_-, \; \pi_{-1})  \lrw (BvB/B,\; \pi_1)
\]
is an anti-Poisson isomorphism. 
\ele

\begin{proof}
It is proved in \cite[Appendix A]{Elek-Lu:BS} that if
$\rho_\sY: (X, \piX) \to (Y, \piY)$ and $\rho_{\sZ}: (X, \piX) \to (Z, \piZ)$ form a weak Poisson pair
and if $X'$ is a coisotropic submanifold of
$(X, \piX)$ such that $\rho_\sY|_{\sX'}: X' \to Y$ is a diffeomorphism, then
$\Phi = \rho_\sZ \circ (\rho_\sY|_{\sX'})^{-1}:  (Y, \piY) \to (Z, \piZ)$
is an anti-Poisson map. Applying the above statement to the weak Poisson pair $(\varpi_-|_{\sG^{v, v}}, \, \varpi|_{\sG^{v, v}})$
in \leref{le-weak-pair} and the coisotropic submanifold $C_\bv$ of $(\Gvv, \pist)$, one 
proves \leref{le-Phi-bv}.
\end{proof}

\bre{re-Phi-uv-v}
With $\Phi_\bv$ defined in \eqref{eq-Phi-bv}, for $u \in W$, let
\begin{equation}\lb{eq-varpi-uv-v}
\varpi^{u, v}_\bv = \Phi_\bv \circ (\varpi_-|_{\sG^{u, v}}): \;\;
\Guv \lrw BvB/B, \;\; b_-c \longmapsto c_\cdot B,
\hs b_- \in B_-,  \, c \in C_\bv.
\end{equation}
It follows from \leref{le-Phi-bv} that 
$\varpi_\bv^{u, v}: (\Guv, \,\pist) \to (BvB/B, \,\pi_1)$
is anti-Poisson.
Consequently, by \leref{le-weak-pair}, one has the weak Poisson pairs
\begin{align*}
&\varpi|_{\sG^{u, v}}: \;\; (\Guv, \,\pist) \lrw (BuB/B, \;\pi_1) \hs \mbox{and} \hs
\varpi_\bv^{u, v}: \;\; (\Guv, \,\pist) \lrw (BvB/B, \;-\pi_1),\\
&\varpi|_{\sSigma^{u, v}}: \;\; (\Sigma^{u, v}, \,\pist) \lrw (BuB/B, \;\pi_1) \hs \mbox{and} \hs
\varpi_\bv^{u, v}: \;\; (\Sigma^{u, v}, \,\pist) \lrw (BvB/B, \;-\pi_1),
\end{align*}
where $\Sigma^{u, v}$ is any symplectic leaf of $\pist$ in $\Guv$.
Note that when $u = v$, $\varpi|_{\sG^{v, v}} = \theta_\bv$ and 
$\varpi^{v, v}_\bv = \tau_\bv$, the source and target maps of the Poisson
groupoid $(G^{\bv, \bv}, \pist)$ over $(BvB/B, \pi_1)$.
\hfill $\diamond$
\ere




\subsection{Commuting Poisson actions of $(G^{\bu, \bu}, \pist)$ and $(G^{\bv, \bv}, \pist)$
on $(\Guv, \pist)$}\lb{subsec-Poisson-bi}
Recall that if $(\G \rightrightarrows Y, \pi_{\ssG})$ is a Poisson groupoid over a Poisson manifold $(Y, \piY)$
with target map $\tau: \G \to Y$, 
a {\it left Poisson action} of $(\G, \pi_{\ssG})$ on a Poisson manifold $(X, \piX)$ is a left 
Lie groupoid $\G$-action on $X$ with a moment map $\nu: X \to Y$ and an action map 
\[
{\bf a}: \;\; \G \ast X \, \stackrel{{\rm def}}{=}\, \{(\gamma, x) \in \G \times X: \, \tau(\gamma) = \nu(x)\} \lrw X
\]
such that ${\rm Graph}({\bf a}) \stackrel{{\rm def}}{=} \{(\gamma, x, {\bf a}(\gamma, x)): (\gamma, x) \in \G \ast X\}$
is a coisotropic submanifold of the Poisson manifold 
$(\G \times X \times X, \, \pi_{\ssG} \times \piX \times (-\piX))$.
In such a case, 
the moment map $\nu: (X, \piX) \to (Y, \piY)$ is automatically Poisson \cite{Liu-Wei-Xu:dirac}.
Note that the moment map $\nu$ is required to be a submersion to ensure that $\G \ast X$ is 
a smooth submanifold of $\G \times X$. Right Poisson actions of Poisson groupoids are similarly defined, where
the moment maps are necessarily anti-Poisson.


Let now $u, v \in W$ and let $\bu, \bv$ be any respective representatives of $u$ and $v$ in $N_\sG(T)$.
Then it is straightforward to check that the groupoid 
$G^{\bu, \bu}$ acts on $G^{u, v}$ on the left with the moment map
$\varpi|_{\sG^{u, v}}: G^{u, v} \to BuB/B$, where 
the action of $g \in G^{\bu, \bu}$ on $x \in G^{u, v}$ with $\tau_{\bu}(g) = \varpi(x)$ is the element
$g\triangleright x \in \Guv$ given by
\begin{equation}\lb{eq-action-1}
g\triangleright x \, \stackrel{{\rm def}}{=}\, cbb' = b_-b_-^\prime c'', \hs\; \mbox{if} \;\;
g = cb = b_-c', \;\; 
x = c'b' = b_-^\prime c^{\prime\prime},
\end{equation}
with
$c, \,c' \in C_\bu, \, c^{\prime\prime} \in C_{\bv}, \,b, \,b' \in B$ and $b_-, \,b_-^\prime \in B_-$. Similarly
the groupoid $G^{\bv, \bv}$ acts on $G^{u, v}$ on the right with the moment map
$\varpi_\bv^{u, v}: G^{u, v} \to BvB/B$ (see \eqref{eq-varpi-uv-v}), and 
the action of $h\in G^{\bv, \bv}$ on $x \in G^{u, v}$ with $\varpi_\bv^{u, v}(x) = \theta_\bv(h)$ is the element
$x \triangleleft h \in \Guv$ given by 
\begin{equation}\lb{eq-action-2}
x \triangleleft h \, \stackrel{{\rm def}}{=}\, c'b'b^{\prime\prime} = b_-^\prime
b_-^{\prime\prime} c^{\prime\prime\prime}, \hs \mbox{if}\; x = c'b' = b_-^\prime c^{\prime\prime},
\;\;
\mbox{and} \;\; h = c^{\prime\prime}b^{\prime\prime} = b_-^{\prime\prime} c^{\prime\prime\prime},
\end{equation}
with
$c' \in C_\bu$, $c^{\prime\prime}, \, c^{\prime\prime\prime} \in C_{\bv}, \,b', \,b^{\prime\prime} \in B$ and $b_-^{\prime}, \,b_-^{\prime\prime} \in B_-$.
One can also check directly that the two groupoid actions commute.

\bthm{thm-Poi-bi}
For any $u, v \in W$ and respective representatives $\bu, \bv \in N_\sG(T)$, 
\eqref{eq-action-1} and \eqref{eq-action-2} are respectively left and right Poisson actions
of the Poisson groupoids $(G^{\bu, \bu}, \pist)$ and $(G^{\bv, \bv}, \pist)$ on $(\Guv, \pist)$.
\ethm

\begin{proof} Consider first the right action of
$(G^{\bv, \bv}, \pist)$ on $(\Guv, \pist)$. Under the Poisson embedding $I_\bv: (B_- v B_-, \pist) \to
((G/B) \times B_-, \pi_1)$, one has $I_\bv(G^{u, v}) = F^{u, v}$ and $I_\bv(\Gvv) = F^{v, v}$ (see
\eqref{eq-Fuv}), and the right action of $G^{\bv, \bv}$ on $\Guv$ corresponds to the 
right action of $F^{v, v}$ on
$F^{u, v}$ by restricting the right Poisson action of the Poisson groupoid
$((G/B) \times B_-, \pi)$ on itself by right multiplication with the target map $\tau$ as the moment map. As
$F^{v, v}$ and $F^{u, v}$ are both Poisson submanifolds of $(G/B) \times B_-$ with respect to $\pi$,
the right action of $(G^{\bv, \bv}, \pist)$ on $(\Guv, \pist)$ is Poisson. 

By replacing $(u, v)$ by $(v, u)$ in the above arguments, one has a right Poisson action
\begin{equation}\lb{eq-vu-uu}
G^{v, u} \times G^{\bu, \bu} \ni (x, \, g) \longmapsto x \triangleleft g \in G^{v, u}, \hs \mbox{if} \;\;\;
\varpi^{v, u}_\bu(x) = \theta_\bu(g).
\end{equation}
One now checks directly that under the Poisson isomorphisms
\[
(\iota^{\bu, \bv})^{-1}: \;\; (G^{v, u}, \, \pist) \lrw (\Guv, \, -\pist) \hs \mbox{and} \hs
\iota_{\bu}: \;\; (G^{\bu, \bu}, \, \pist) \lrw (G^{\bu, \bu}, \, -\pist),
\]
where the Poisson isomorphism $\iota^{\bu, \bv}: (\Guv, \, \pist) \to (G^{v, u}, \, -\pist)$ is given in 
\eqref{eq-iota-uv},  the right groupoid action of $G^{\bu, \bu}$ on $G^{v, u}$ in \eqref{eq-vu-uu}
becomes precisely the left action of the groupoid $G^{\bu, \bu}$ on $G^{u, v}$ given in \eqref{eq-action-1}.
This shows that the groupoid action in \eqref{eq-action-1} is Poisson.
\end{proof}

\sectionnew{Symplectic groupoids associated to double Bruhat cells}\lb{sec-symplectic}

\subsection{Symplectic leaves in $\Guv$}\lb{subsec-leaves-Guv}
To describe the symplectic leaves of $\pist$ in $G$,
it is enough to describe the symplectic leaves in the double Bruhat cells, as the 
latter are the $T$-orbits of symplectic leaves
of $\pist$ in $G$. For $u, v \in W$, and for any symplectic leaf $\Sigma$ of $\pist$ in $\Guv$, let
$T_{\sSigma} = \{t \in T: \Sigma \, t = \Sigma\}$. As $T$ acts transitively on the set of all symplectic leaves
of $\pist$ in $\Guv$, $T_\sSigma$ is independent of $\Sigma \subset \Guv$. We define the 
{\it leaf-stabilizer of $T$ in $\Guv$} to be
\begin{equation}\lb{eq-T-stab}
T^{u, v}_{\rm stab} = T_\sSigma, 
\end{equation}
where $\Sigma$ is any symplectic leaf of $\pist$ in $\Guv$.
When $G$ is simply connected, symplectic leaves of $\pist$ in each $\Guv$ 
are determined by Kogan and Zelevinsky in \cite{k-z:leaves} using specially chosen representatives 
in $N_\sG(T)$ of elements in $W$. In this section, for $G$ simply connected,
we adapt the results in \cite{k-z:leaves} to
describe the symplectic leaves of $\pist$ in $G$ using arbitrary choices of 
representatives of elements in $W$, and we describe the leaf-stabilizers of $T$ in the double Bruhat cells. We
also extend some results from \cite{k-z:leaves} to the case when $G$ is 
not necessarily simply connected. 

Assume first that $G$ is connected but not necessarily simply connected. The action of the Weyl group on $T$ will be denoted as $t^v = \bv^{-1} t \bv$, where $v \in W, t \in T$, and $\bv$ is any representative of $v$ in $N_\sG(T)$. For $u, v \in W$, let 
\[
T^{u, v} = \{(t^u)^{-1} t^v: t \in T\}.
\] 
Fix $u, v \in W$ and let $\bu, \bv$ be any representatives of $u$ and $v$ in $N_\sG(T)$, respectively. Note that 
\[
\bu^{-1}BuB = \bu^{-1} C_\bu B \subset N_-TN \hs \mbox{and} \hs B_-vB_- \bv^{-1} = B_- C_\bv \bv^{-1}\subset N_-TN,
\]
and recall that for $g \in N_-TN$, we write $g = [g]_-[g]_0[g]_+$, where $[g]_- \in N_-, \, [g]_0 \in T, \, [g]_+ \in N$. For $t \in T$, define 
\begin{equation}\lb{eq-Suvt}
S^{\bu, \bv}_{[t]} =\left\{g \in \Guv: \;
\left[\bu^{-1}g\right]_0 \left[g \, \bv^{-1}\right]_0^v \in tT^{u, v}\right\},
\end{equation}
where $[t]$ denotes the image of $t$ in $T/T^{u, v}$. Define the map 
\begin{equation}\lb{eq-chi-uv}
\chi: \;\; \Guv \lrw T/T^{u, v}, \;\;\; \chi(g) = [\bu^{-1} g]_0 [g\, \bv^{-1}]_0^v\,  T^{u, v} \in T/T^{u, v}, \hs \hs g \in \Guv.
\end{equation}
Then clearly $\Suvt = \chi^{-1}([t])$ for $t \in T$, a level set of $\chi$. One also has
\begin{equation}\lb{eq-chi-a}
\chi(ga) = [a]^2 \chi(g), \hs g \in \Guv, \; a \in T.
\end{equation}
The following \leref{le-Suv} is proved in \cite[Proposition 3.1]{k-z:leaves} (neither the assumption that $G$ be simply-connected nor the the special way  of choosing representatives of Weyl group elements in $N_\sG(T)$ made in \cite{k-z:leaves} is needed in its proof).

\ble{le-Suv} \cite[Proposition 3.1]{k-z:leaves}
The symplectic leaves of $\pist$ in $\Guv$ are the connected components of the sets $S^{\bu, \bv}_{[t]}$, $t \in T$. Moreover, for any $t_1, t_2, t \in T$, $S_{[t_1]}^{\bu, \bv} = S_{[t_2]}^{\bu, \bv}$ if and only if $[t_1] = [t_2]$, and $S_{[t_1]}^{\bu, \bv} t = S_{[t_1t^2]}^{\bu, \bv}$.
\ele


Assume now that $G$ is simply-connected, and let $\Gamma \subset \Delta_+$ be the set of 
simple roots. For $\alpha \in \Gamma$, let $\omega_\alpha \in \Hom(T, \Cset^\times)$ be the 
corresponding fundamental weight, and let $\Delta_\al$ be the corresponding generalized principal minor
\cite{Fomin-Zelevinsky:double, k-z:leaves}, which is a regular function on $G$ whose restriction to 
$N_-TN$ is given by $\Delta_\al(g) = [g]_0^{\omega_\al}$.
For $u, v \in W$, let $I(u, v) = I(u) \cap I(v)$, where 
\[
I(u) = \{\alpha \in \Gamma: \; u(\omega_\al) = \omega_\al\} =\Gamma\backslash\{\alpha_1, \ldots, \al_l\}
\]
for any reduced word $u = s_{\alpha_1} s_{\alpha_2} \cdots s_{\alpha_l}$, and
define the maps $\delta, \, \delta^2: G \to \Cset^{|I(u, v)|}$ by
\[
\delta(g) = \{\Delta_\al(g): \, \alpha \in I(u, v)\}
\hs \mbox{and} \hs \delta^2(g) = \{(\Delta_\al(g))^2: \, \alpha \in I(u, v)\}.
\]
We now modify the results from \cite{k-z:leaves} to give a description of the
connected components of $S^{\bu, \bv}_{[t]}$,
and thus also of the symplectic leaves of $\pist$ in $\Guv$.

\bpr{pr-leaf-uv} 
Assume that $G$ is simply connected. Let $u, v \in W$ and let  $\bu$ and $\bv$ be any
respective representatives of $u$ and $v$ in $N_\sG(T)$. Then for any $t \in T$, the restriction of
$\delta^2$ to $\Suvt$ is a constant map, or, more precisely,
\begin{equation}\lb{eq-delta-g}
(\Delta_\al(g))^2 = \Delta_\al(\bu) \Delta_\al(\bv) \,t^{\omega_\al}, \hs \;\; \forall \;\; g \in \Suvt.
\end{equation}
The connected components of $\Suvt$  
are the $2^{|I(u, v)|}$ (all of which non-empty) level sets of the map $\delta: \Suvt \to (\Cset^\times)^{|I(u, v)|}$.
\epr

\begin{proof}
By first choosing a set $\{e_\alpha \in \g_\alpha, \, f_\alpha \in \g_{-\alpha}, \, \alpha^\vee
\in \h: \alpha \in \Gamma\}$ of Chevalley generators of $\g$
which determines Lie group homomorphisms $\phi_\alpha: SL(2, \Cset) \to G$ for each $\alpha \in \Gamma$, one can
choose the representative $\tilde{s}_\alpha$ of $s_\alpha$ in $N_\sG(T)$ to be
$\tilde{s}_\alpha = \phi_\alpha \left(\begin{array}{cc} 0 & -1 \\ 1 & 0\end{array}\right)$ for each $\alpha \in
\Gamma$. For $w \in W$ and any reduced word $w = s_{\alpha_1} s_{\alpha_2}\cdots s_{\alpha_l}$ of $w$, the element
$\tilde{w} = \tilde{s}_{\alpha_1} \tilde{s}_{\alpha_2} \cdots \tilde{s}_{\alpha_l}$ is then a representative of $w$ in $N_\sG(T)$ 
independent of the choice of the reduced word. Moreover \cite[Lemma 6.1]{M-R:parametrizations},
$\Delta_\al(\tilde{w})=1$ if $\al \in I(w)$.  Define 
\[
S^{u, v}_{e} = \left\{g \in \Guv: \; \left[\tilde{u}^{-1}g\right]_0 \left[g \,
\widetilde{v^{-1}}\right]_0^v \in T^{u, v}\right\}.
\]
By \cite[Theorem 2.3, Corollary 2.5, Lemma 3.2]{k-z:leaves}, 
$[\tilde{u}^{-1} g]_0^{\omega_\al} = \pm 1$ for all $g \in S^{u, v}_{e}$ and $\al \in I(u, v)$, 
and $S^{u, v}_{e}$ has $2^{|I(u, v)|}$ connected components
$S^{u, v}_{e} =\bigsqcup_{\epsilon} S^{u, v}_{{e}}(\epsilon)$, where $\epsilon$ runs over 
the set of all sign functions
$\epsilon: I(u, v) \to \{\pm 1\}$ on $I(u, v)$, and 
\[
S^{u, v}_e(\epsilon) = \left\{g \in S^{u, v}_e: \; 
\;\, \left[\tilde{u}^{-1}g\right]_0^{\omega_\al} = \epsilon(\alpha), \; \forall \, \alpha \in I(u,v)\right\}.
\]
Let $t_0, t_1 \in T$ be such that $\tilde{u}=t_0\bu$ and $\widetilde{v^{-1}} \bv = t_1$. One checks directly that
for $g \in \Guv$,
\[
\left[\tilde{u}^{-1}g\right]_0 = \left[\bu^{-1}g\right]_0 (t_0^{-1})^u \hs \mbox{and} \hs
\left[g \,\widetilde{v^{-1}}\right]_0^v = \left[g \,\bv^{-1}\right]_0^vt_1.
\]
It follows that for any $t \in T$ and $a \in T$ with $a^2 = (t_0^{-1})^u t_1t$, one has
\[
\Suvt = \left\{g \in \Guv: \; \left[\tilde{u}^{-1}g\right]_0 \left[g \,
\widetilde{v^{-1}}\right]_0^v \in (t_0^{-1})^u  t_1 tT^{u, v}\right\} = S^{u, v}_e a.
\]
Note now (see \cite[(3.12)]{k-z:leaves}) that $\left[\tilde{u}^{-1}g\right]_0^{\omega_\al} = \Delta_\al(g)$
for any $\al \in I(u)$ and  $g \in BuB$. It follows that  $\Delta_\al(g) = \pm 1$ for all 
$g \in S^{u, v}_e$ and $\al \in I(u, v)$. Consequently, for all $g \in \Suvt$ and $\al \in I(u, v)$, 
\[
(\Delta_\al(g))^2 = \Delta_\al(a^2) = t_0^{-\omega_\al} t_1^{\omega_\al} t^{\omega_\al} = \Delta_\al(\bu) 
\Delta_\al(\bv) t^{\omega_\al}.
\]
As there is one connected component of $S^{u, v}_e$ for each sign function $\epsilon$ on $I(u, v)$,
the connected  components of $\Suvt =S^{u, v}_e a$ are precisely the
$2^{|I(u, v)|}$ level sets, all of which non-empty, of the map $\delta: \Suvt \to (\Cset^\times)^{|I(u, v)|}$.
\end{proof}

Recall the map $\chi: \Guv \to T/T^{u, v}$ defined in \eqref{eq-chi-uv}. The following \coref{co-leaf-any} is also proved in 
\cite[Corolary 4.5]{Yakimov:Spectra}.

\bco{co-leaf-any}
For 
any $g_0 \in \Guv$, the symplectic leaf
$\Sigma^{g_0}$ of $\pist$ through $g_0$ is given by
\begin{equation}\lb{eq-Sigma-g0}
\Sigma^{g_0} = \{g \in \Guv: \; \chi(g) = \chi(g_0) \;\, \mbox{and} \;\, \delta(g) = \delta(g_0)\}.
\end{equation}
\eco

\bre{re-Delta-g}
Using the decompositions $BuB = C_\bu B$ and $B_-vB_- = B_- C_\bv$, one can describe the maps $\chi$ and $\delta$ on $\Guv$
more explicitly. Indeed, writing an element $g \in \Guv$ as $g = c t n = n_-t_-c'$, where $c \in C_\bu$, $c' \in C_{\bv}$, 
$t, t_- \in T$, $n \in N$ and $n_- \in N_-$, one has $\chi(g) = [tt_-^v] \in T/T^{u, v}$ and 
$\Delta_\al(g) = t^{\omega_\al} \Delta_\alpha(\bu)$ for all $\al \in I(u)$. 
\hfill $\diamond$
\ere

When $u = v$, one has $T^{u, v} = \{e\}$. As a special case of \coref{co-leaf-any}, one has

\bco{co-leaf-bv} 
Assume that $G$ is simply connected. Let $v \in W$ and let $\bv$ be any
representative of $v$ in $N_\sG(T)$. 
Then the symplectic leaf $\Sigma^\bv$ of $\pist$ in $G$ through $\bv$ is given by
\begin{align}\lb{eq-Sigma-bv}
\Sigma^\bv &= \left\{g \in \Gvv: \;
\left[\bv^{-1}g\right]_0 \left(\left[g \, \bv^{-1}\right]_0\right)^v =e, 
\; \Delta_\al(g) = \Delta_\al(\bv)\,\; \forall 
\alpha \in I(v)\right\}\\
\lb{eq-Sigma-bv-1}
&= \left\{g \in \Gvv: \;
\left[\bv^{-1}g\right]_0 \left(\left[g \, \bv^{-1}\right]_0\right)^v =e, 
\; \left[\bv^{-1}g\right]_0^{\omega_\al} = 1\,\; \forall 
\alpha \in I(v)\right\}
\end{align}
\eco

Still assuming that $G$ is simply connected, let
\[
\tilde{T}^{u, v} = \{t \in T: \;\, t^{\omega_\al} = 1 \;\forall \; \al \in I(u, v)\}.
\]
 It is clear that
$T^{u, v} \subset \tilde{T}^{u, v}$. As a direct consequence of \coref{co-leaf-any} and \eqref{eq-chi-a}, one has 

\bco{co-leaf-stabilizer}
Assume that $G$ is simply connected. Then for $u, v \in W$, the leaf-stabilizer of $T$ in $\Guv$ is given by
$T_{\rm stab}^{u, v}=\{t \in \tilde{T}^{u, v}: \; t^2 \in T^{u, v}\}$. 
\eco

Returning now to the connected semisimple complex Lie group $G$ which may not be simply connected, 
let $\hat{G}$ be the connected and simply connected cover of $G$, and let $\kappa: \hat{G} \to G$ be the covering 
map with $\ker \kappa = Z$, a subgroup of the center of $\hat{G}$. Denoting by $\hat{\pi}_{\rm st}$ the multiplicative Poisson structure on $\hat{G}$ 
defined by the same $r$-matrix $r_{\rm st} \in \gotg$, the map
$\kappa: (\hat{G}, \hat{\pi}_{\rm st})\to (G, \pist)$ is then Poisson.
For $\hat{g} \in \hat{G}$ and $g \in G$, let again $\Sigma^{\hat{g}} \subset \hat{G}$ and $\Sigma^{g} \subset G$
respectively denote the symplectic leaves of $\hat{\pi}_{\rm st}$ and $\pist$ through $\hat{g}$ and $g$. 
Let $\hat{T} = \kappa^{-1}(T)$, a maximal torus of $\hat{G}$. By \coref{co-leaf-stabilizer},
the leaf-stabilizer of $\hat{T}$ in $\hat{G}^{u, v}$ is
\[
\hat{T}_{\rm stab}^{u, v} =\{\hat{a} \in \hat{T}: \; \hat{a}^{\omega_\al} = 1\;\forall \;
\al \in I(u, v)\;\; \mbox{and} \;\;  \hat{a}^2 \in 
\hat{T}^{u, v}\}.
\]
Let $Z_{u, v} = Z \cap \hat{T}_{\rm stab}^{u, v} = \{z \in Z:\; z^{\omega_\al} = 1\;\forall \;
\al \in I(u, v)\;\; \mbox{and} \;\;  z^2 \in 
\hat{T}^{u, v}\}.$

\ble{le-tv-bv} 
For any $\hat{g} \in \hat{G}^{u, v}$, one has $\kappa\left(\Sigma^{\hat{g}}\right) = \Sigma^g$, where $g = \kappa(\hat{g}) \in \Guv$, and $\kappa: \Sigma^{\hat{g}} \to \Sigma^g$ is the quotient map $\Sigma^{\hat{g}} \to \Sigma^{\hat{g}}/Z_{u, v}$, where $Z_{u, v}$ acts on $\Sigma^{\hat{g}}$ by multiplication.
\ele

\begin{proof}
As $\kappa: (\hat{G}, \hat{\pi}_{\rm st})\to (G, \pist)$ is a local Poisson diffeomorphism, and as
$\Sigma^{\hat{g}}$ is connected, we have $\kappa(\Sigma^{\hat{g}}) \subset \Sigma^g$.
To show that $\Sigma^g \subset \kappa(\Sigma^{\hat{g}})$, let $h \in \Sigma^g$ and
let $\gamma: [0, 1] \to \Sigma^g$ be any smooth path in $\Sigma^g$ such that $\gamma(0) = g$ and $\gamma(1) = h$.
Let $\hat{\gamma}: [0, 1] \to \hat{G}$ be the unique lifting of $\gamma$ such that $\hat{\gamma}(0) = \hat{g}$. 
Again as $\kappa$ is a local Poisson diffeomorphism, $\hat{\gamma}$ is tangent to the symplectic leaf
through $\hat{\gamma}(x)$ for every $x \in [0, 1]$. Thus 
$\hat{\gamma}([0, 1]) \subset \Sigma^{\hat{g}}$. This shows that
$\kappa\left(\Sigma^{\hat{g}}\right) = \Sigma^g$.

Clearly the $Z_{u, v}$-orbits in $\Sigma^{\hat{g}}$ are contained in the fibers of 
$\kappa: \Sigma^{\hat{g}} \to \Sigma^g$. Suppose that
$\hat{h}, \hat{k} \in \Sigma^{\hat{g}}$ are in the same fiber of $\kappa: \Sigma^{\hat{g}} \to \Sigma^g$. Then 
$\hat{h}z = \hat{k}$ for some $z \in Z$. As $\Sigma^{\hat{g}}z$ and $\Sigma^{\hat{g}}$ are both symplectic leaves
of $\hat{\pi}_{\rm st}$ and have now a non-empty intersection, $\Sigma^{\hat{g}}z =\Sigma^{\hat{g}}$, and
thus $z \in Z_{u, v}$. 
\end{proof}

\bre{re-convering}
Same arguments as in the proof of \leref{le-tv-bv} show that if $\kappa: (X, \piX) \to 
(Y, \piY)$ is a covering map that is also Poisson, then the images under $\kappa$ of the symplectic leaves 
of $(X, \piX)$ are precisely all the symplectic leaves of $(Y, \piY)$.
\hfill $\diamond$
\ere

\ble{le-leaf-stab} 
For any $u, v \in W$, the leaf-stabilizer of $T$ in $\Guv$ is given by $T_{\rm stab}^{u, v} = \kappa\left(\hat{T}_{\rm stab}^{u, v}\right)$.
\ele

\begin{proof} Let $\hat{\Sigma}$ be a symplectic of $\hat{\pi}_{\rm st}$ in $\hat{G}^{u, v}$, and let
$\Sigma =\kappa(\hat{\Sigma})$. 
If $\hat{a} \in \hat{T}_{\rm stab}^{u, v}$, then it follows from  
$\hat{\Sigma} \hat{a} = \hat{\Sigma}$ that $\Sigma \kappa(\hat{a}) =\Sigma$, so $\kappa(\hat{a}) \in T^{u, v}_{\rm stab}$. Conversely,
let $a \in T^{u, v}_{\rm stab}$
and choose any $\hat{a} \in \kappa^{-1}(a)$. Let $\hat{g} \in \hat{\Sigma}$. Then $\kappa(\hat{g}\hat{a}) 
\in\Sigma a = \Sigma$, so $\kappa(\hat{g}\hat{a}) = \kappa(\hat{g}^\prime)$ for some 
$\hat{g}^\prime \in \hat{\Sigma}$. Let $z \in Z$ be such that $\hat{g}\hat{a} z = \hat{g}^\prime$.
As  $\hat{\Sigma} \hat{a}z$ and $\hat{\Sigma}$ are two symplectic leaves of $\hat{\pi}_{\rm st}$
and have a non-empty intersection, one must have $\hat{\Sigma} \hat{a}z =\hat{\Sigma}$, and thus
$a = \kappa(\hat{a}z) \in \kappa\left(\hat{T}_{\rm stab}^{u, v}\right)$.
\end{proof}

Recall from \leref{le-Suv} that symplectic leaves of $\pist$ in $\Guv$ are the connected components of
the sets $\Suvt$ given in \eqref{eq-Suvt}, where $t \in T$. Define
\[
T^{(2)} = \{a \in T: \; a^2 = e\}.
\]
It is clear that for each $t \in T$, $\Suvt$ is invariant under left translation by elements in $T^{(2)}$.

\ble{le-T2-transtitive}
For any $t \in T$, the induced action of $T^{(2)}$ on the set of all symplectic leaves of $\pist$ in $\Suvt$ is transitive.
\ele

\begin{proof}
Let $\Sigma$ and $\Sigma^\prime$ be any two symplectic leaves of $\pist$ in $\Suvt$, and let 
$\hat{\Sigma}$ and $\hat{\Sigma}^\prime$ be two symplectic leaves of $\hat{\pi}_{\rm st}$ in $\hat{G}^{u, v}$
such that $\kappa(\hat{\Sigma}) = \Sigma$ and $\kappa(\hat{\Sigma}^\prime) = \Sigma^\prime$. 
Let
$[\kappa]: \hat{T}/\hat{T}^{u, v} \to T/T^{u, v}$ 
be the group homomorphism induced by $\kappa: \hat{T} \to T$. Then the fibers of $[\kappa]$ are
the $Z$-orbits in $\hat{T}/\hat{T}^{u, v}$ by multiplication.
Let $\hat{u}$ and $\hat{v}$ be any 
respective representatives of $u$ and $v$ in 
$N_{\scriptscriptstyle{\hat{G}}}(\hat{T}) \subset \hat{G}$. Recalling the map 
$\hat{\chi}: \hat{G}^{\hat{u}, \hat{v}} \to \hat{T}/\hat{T}^{u, v}$ defined as in \eqref{eq-chi-uv}, one has
\[
[\kappa](\hat{\chi}(\hat{\Sigma})) = [\kappa](\hat{\chi}(\hat{\Sigma}^\prime)) = [t].
\]
Thus there exists $z \in Z$ such that
$z\hat{\chi}(\hat{\Sigma}) = \hat{\chi}(\hat{\Sigma}^\prime)$. Let $\hat{a} \in \hat{T}$ be such that
$\hat{a}^2 = z$. Then $\hat{\chi}(\hat{\Sigma}\hat{a}) = \hat{\chi}(\hat{\Sigma}^\prime)$.  
By \prref{pr-leaf-uv} (see also the proof of \cite[Theorem 2.3]{k-z:leaves}), the group
$\hat{T}^{(2)} = \{\hat{x} \in \hat{T}: \, \hat{x}^2 = e\}$ acts transitively on the set of the symplectic leaves of
$\hat{\pi}_{\rm st}$ in any level set of $\hat{\chi}$. Thus there exists $\hat{x} \in \hat{T}^{(2)}$ such that
$\hat{\Sigma}\hat{a}\hat{x} = \hat{{\Sigma}}^\prime$. Let $a =\kappa(\hat{a} \hat{x}) \in T$. Then $a \in T^{(2)}$
and $\Sigma a = \Sigma^\prime$. 
\end{proof}

\bre{re-number}
It follows from \leref{le-T2-transtitive} and \leref{le-leaf-stab} that
for $t \in T$, the number of symplectic leaves of $\pist$ in $\Suvt$ is equal to
$\displaystyle \left| T^{(2)}/T^{(2)} \cap T_{\rm stab}^{u, v}\right|$. As $T^{(2)}$ is a $2$-group,
the number of symplectic leaves of $\pist$ in $\Suvt$ is always a power of $2$.
\hfill $\diamond$
\ere

\subsection{The symplectic leaf $\Sigma^\bv$ as a symplectic groupoid}\lb{subsec-leaf-groupoid}
Let now $(G, \pist)$ be any standard complex semisimple Poisson Lie group, where $G$ is connected but
not necessarily simply connected.
Let $v, u \in W$ and let $\bu$ and $\bv$ be any respective representatives of $u$ and $v$ in $N_\sG(T)$. 
One then has the Poisson groupoid $(G^{\bu, \bu}, \pist)$ over $(BuB/B, \pi_1)$ and the Poisson groupoid
$(G^{\bv, \bv}, \pist)$ over $(BvB/B, \pi_1)$. Recall their commuting (left and right) Poisson actions  
on $(\Guv, \pist)$, respectively given in \eqref{eq-action-1} and \eqref{eq-action-2}. 

\bthm{symp-leaves} 
1) The symplectic leaf $\Sigma^\bv$ of $\pist$ through $\bv$ is a Lie subgroupoid of $G^{\bv, \bv}$.
Consequently, $(\Sigma^\bv, \pist)$ is a symplectic groupoid over $(BvB/B, \pi_1)$;

2) For any symplectic leaf $\Sigma^{u, v}$ of $\pist$ in $\Guv$, the two commuting Poisson actions 
in \eqref{eq-action-1} and \eqref{eq-action-2} restrict to Poisson actions of the 
symplectic groupoids $(\Sigma^\bu,  \pist)$ and $(\Sigma^\bv,  \pist)$
on the symplectic manifold $(\Sigma^{u, v}, \pist)$. 
\ethm

\begin{proof}  
Assume first that $G$ is simply connected. Consider the action in \eqref{eq-action-1}. 
Assume that $g \in \Sigma^\bu$ and $x \in \Sigma^{u, v}$
be such that $\tau_\bu(g) = \varpi(x)$, and write 
$g = c t n = n_- t_- c'$ and $x = c' t' n' = n_-^\prime t_-^\prime c^{\prime\prime}$,
where $c, c' \in C_\bu, c^{\prime\prime} \in C_\bv$, $t, t_-, t', t_-^\prime \in T$, $n,n' \in N$, and
$n_-, n_-^\prime \in N_-$. Then
$g \triangleright x = ctnt'n' = n_-t_-n_-^\prime t_-^\prime c^{\prime\prime}$. 
By \prref{pr-leaf-uv}, $t t_-^u = e$, $t^{\omega_\al} = 1$ for all $\al \in I(u, v)$, and 
$\Sigma^{u, v} = \{h \in \Guv: \chi(h) = \chi(x), \; \Delta_\al(h) = \Delta_\al(x) \; \forall \; \al \in I(u, v)\}.$
By the definitions of the map  $\chi$ and the functions $\Delta_\al$ (see \reref{re-Delta-g}),
\[
\chi(g \triangleright x) = [tt' (t_-t_-^\prime)^v] = [tt_-^u t' (t_-^\prime)^v (t_-^u)^{-1}t_-^v] = 
[t' (t_-^\prime)^v] = \chi(x) \in T/T^{u, v},
\]
and for every $\al \in I(u, v)$, $\Delta_\al(g \triangleright x) = (tt')^{\omega_\al} \Delta_\al(\bu) = (t')^{\omega_\al} \Delta_\al(\bu) = \Delta_\al(x)$. 
Thus $g \triangleright x \in \Sigma^{u, v}$. 
Similarly, one shows that for all 
$x \in \Sigma^{u, v}$ and $h \in \Sigma^\bv$ with $\varpi_\bv^{u, v}(x) = \theta_\bv(h)$ one has
$x \triangleleft h \in \Sigma^{u, v}$.  
Applying to the special case of $u = v, \bu = \bv$ and $\Sigma^{u, v} = \Sigma^\bv$,
it shows in particular that $\Sigma^\bv$ is closed under
the groupoid multiplication of $G^{\bv, \bv}$. It is easy to see that $\Sigma^\bv$ is closed under 
the groupoid inverse of $G^{\bv, \bv}$. By \leref{le-Sigma-uv}, both 
$\varpi|_{\sSigma^{u, v}}: \Sigma^{u, v} \to BuB/B$ and 
$\varpi_\bv^{u, v}|_{\sSigma^{u, v}}: \Sigma^{u, v} \to BvB/B$ are submersions. 
Thus $\Sigma^\bv$ is a Lie subgroupoid of $G^{\bv, \bv}$  and the two actions
in \eqref{eq-action-1} and \eqref{eq-action-2} restrict to Poisson  actions of  
the symplectic groupoids $(\Sigma^\bu, \pist)$ and $(\Sigma^\bv, \pist)$
on the symplectic manifold $(\Sigma^{u, v}, \pist)$.  

For an arbitrary $G$, let $\hat{G}$ be the simply connected cover of $G$ with
$\kappa: \hat{G} \to G$ the covering map and multiplicative Poisson structure $\hat{\pi}_{\rm st}$,
and choose any $\hat{u}, \hat{v} \in \hat{G}$ such that $\kappa(\hat{u}) = \bu$ and $\kappa(\hat{v}) = \bv$.
Let $Z = \ker \kappa$, and let $\hat{\Sigma}^{u, v}$ be any symplectic leaf of $\hat{\pi}_{\rm st}$ 
such that $\kappa(\hat{\Sigma}^{u, v}) = \Sigma^{u, v}$. By \leref{le-tv-bv},
the symplectic groupoids $(\Sigma^\bu, \pist)$ and
$(\Sigma^\bv, \pist)$ are the respective quotients of the 
symplectic groupoids $(\Sigma^{\hat{u}}, \hat{\pi}_{\rm st})$ and
$(\Sigma^{\hat{v}}, \hat{\pi}_{\rm st})$ by $Z_{u, u}$ and $Z_{v, v}$, and that
$\kappa: \hat{\Sigma}^{u, v} \to \Sigma^{u, v}$ is the quotient map by $Z_{u, v}$. It is easy to 
see that $Z_{u, u} \subset Z_{u, v}$ and $Z_{v, v} \subset Z_{u, v}$. 
Statements 1) and 2) for $G$ now follow from the corresponding statements for $\hat{G}$.
\end{proof}

\bre{re-all-leaves} 
Let  $u, v \in W$, and let $\Sigma^u \subset G^{u, u}, \Sigma^{u, v} \subset \Guv$, and
$\Sigma^v \subset \Gvv$ be arbitrary symplectic leaves of $\pist$. As
$\Sigma^u = \Sigma^\bu$ and $\Sigma^v = \Sigma^\bv$ for some representatives 
of $\bu$ and $\bv$,  we conclude that $\Sigma^u$ and $\Sigma^v$ are symplectic groupoids, respectively over
$(BuB, \pi_1)$ and $(BvB/B, \pi_1)$, acting by commuting Poisson actions from the left and right 
on the symplectic groupoid $(\Sigma^{u, v}, \pist)$.
\hfill $\diamond$
\ere

\bex{sl2-exa}
Let $G = SL(2, \Cset)$, where the pair $(B, B_-)$ consists of the subgroups of respectively upper and lower triangular matrices, and where $\la x_1, x_2 \ra_\g = \tr(x_1x_2)$, $x_1, x_2 \in \sl(2, \Cset)$.  Writing $g \in G$ as $g = \left(\begin{array}{cc}g_{11} & g_{12} \\g_{21} & g_{22}\end{array}\right)$, the Poisson brackets between the coordinate functions are 
\begin{align*}
&\{g_{11}, g_{12}\} = g_{11}g_{12}, \hs \{g_{11}, g_{21}\} = g_{11} g_{21}, \hs \{g_{12}, g_{22}\} = g_{12} g_{22},\\
&\{g_{21}, g_{22}\}=  g_{21} g_{22}, \hs \{g_{11}, g_{22}\} = 2g_{12}g_{21}, \hs \{g_{12}, g_{21}\} = 0.
\end{align*}
Let $\bs = \left(\begin{array}{cc} 0 & -1 \\ 1 & 0\end{array}\right)$, so that $C_{\bs} = \left\{\left(\begin{array}{cc} z & -1 \\ 1 & 0\end{array}\right): z \in {\mathbb C}\right\}$. Then
\[
G^{s, s} = \left\{\left(\begin{array}{cc} az & a^{-1}(abz-1) \\ a & b\end{array}\right): \;a, \,b, \,z \in \Cset, \;
 a \neq 0, \; abz-1 \neq 0\right\},
\]
with the Poisson structure given by
$\{z, a\} = za,\, \{z, b\} = a^{-1}(abz-2), \, \{a, b\} = ab.$
Let $\chi = a^2(1-abz)^{-1}$. The groupoid structure on $G^{\bs, \bs}$ over $\Cset$ is given by 
\begin{align*}
&\text{source map}: \; \theta_\bs(z, a, b) = z,\\
&\text{target map}: \;\tau_\bs(z, a, b) = \chi z,\\
&\text{inverse map}: \; \iota_\bs(z, a, b) = (\chi z, \, a^{-1}, \, -b),\\
&\text{identity bisection}: \; z \mapsto (z, 1, 0), \; z \in \Cset,\\
&\text{multiplication}: \mu_\bs((z_1, a_1, b_1), (z_2, a_2, b_2)) = (z_1,  a_1a_2,  a_1b_2 + b_1a_2^{-1}),
\;\text{if} \; z_2 = \tau_\bs(z_1,a_1, b_1).
\end{align*} 
Note that $\chi$ is a Casimir function on $G^{s, s}$ and the symplectic leaves in $G^{s, s}$ are precisely given by the (non-zero) level sets of $\chi$. Hence the symplectic leaf $\Sigma_\bs$ of $\pist$ through $\bs \in G^{\bs, \bs}$ is
$\displaystyle 
\Sigma_{\bs} = \left\{\left(\begin{array}{cc} az & -a \\ a & b\end{array}\right): \;a, \,b, \,z \in \Cset, \;
a \neq 0, \; a^2 = 1-abz\right\}$. Identify $\Sigma_\bs$ with 
\begin{equation}\lb{Sigma-defn}
\Sigma = \left\{\left(\begin{array}{cc}pt & -t \\t & -qt\end{array}\right): \;\;
(p, q, t) \in \Cset^3, \; t^2(1-pq) = 1 \right\}.
\end{equation}
The induced (non-degenerate) Poisson structure on $\Sigma$ is given by
\begin{equation} \lb{bracket-Sigma}
\{p, q\} = 2(1-pq), \hs \{p, t\} = pt, \hs \{q, t\} = -qt, 
\end{equation}
and the induced symplectic groupoid structure on $\Sigma$ is given by
\begin{align*}
&\mbox{source map}: \; \theta(p, q, t) = p,\\
&\mbox{target map}: \;\tau(p, q, t) = p,\\
&\mbox{inverse map}: \; \iota(p, q, t) = (p, \, -qt^2,\, t^{-1}),\\
&\mbox{identity section}: \;\epsilon(p) = (p, 0, 1),\\
&\mbox{multiplication}: \; \mu((p_1, q_1, t_1), (p_2, q_2, t_2)) = (p_1, q_2+q_1t_2^{-2}, t_1t_2)
\;\mbox{when} \; p_1=p_2.
\end{align*}
Note that $\theta^{-1}(0)=\tau^{-1}(0)$ is isomorphic to the non-connected abelian Lie group $\Cset \times {\mathbb{Z}}_2$.

Consider now the group $PSL(2, \Cset)$, and write its elements as $[g]$, where $g \in SL(2, \Cset)$.
Then the symplectic leaf of $\pist$ through $[\bs] \in PSL(2, \Cset)$ is parametrized by the surface 
\[
\Sigma_0 = \{(p, q) \in \Cset^2: 1-pq \neq 0\} \cong
\left\{\left[\left(\begin{array}{cc}p & -1 \\1 & -q\end{array}\right)\right]: \;\;
(p, q) \in \Cset^2, \; 1-pq \neq 0 \right\},
\]
with the  Poisson structure $\{p, q\} = 2(1-pq)$ and the 
groupoid structure given by 
\begin{align*}
&\text{source map}: \; \theta(p, q) = p,\\
&\text{target map}: \;\tau(p, q) = p,\\
&\text{inverse map}: \; \iota(p, q, t) = (p, \; q(pq-1)^{-1}),\\
&\text{identity bisection}: \; \{(p, 0) : p \in \Cset \},\\
&\text{multiplication}: \; \mu((p_1, q_1), \,(p_2, q_2)) = (p_1, \;q_2+q_1(1-p_2q_2)),
\; \text{if} \; p_1=p_2.
\end{align*}
Note the Lie group isomorphisms
$\theta^{-1}(0)=\tau^{-1}(0) \cong \Cset$ and $\theta^{-1}(p)=\tau^{-1}(p) \cong \Cset^\times$ for $p \neq 0$.
\hfill $\diamond$
\eex

\bex{sl3-exa}
Let $G = SL(3, \Cset)$, with $B$, $B_-$ respectively the subgroups of upper and lower triangular matrices and the
bilinear form $\la x_1, x_2 \ra_\g = \tr(x_1x_2)$ on $\sl(3, \Cset)$. Let $s_1, s_2$ be the
two generators of the Weyl group $W$, identified  with 
the symmetric group $S_3$. Let $v = s_1s_2$, 
$$
\bs_1 = \left(\begin{array}{ccc}0 & -1 & 0 \\1 & 0 & 0 \\0 & 0 & 1\end{array}\right), \hs \bs_2 = \left(\begin{array}{ccc}1 & 0 & 0 \\0 & 0 & -1 \\0 & 1 & 0\end{array}\right), \hs \text{and} \hs \bv = \bs_1 \bs_2. 
$$
Let $\Sigma_{\bs_1}, \Sigma_{\bs_2}, \Sigma_\bu$ be the symplectic leaves of $\pist$ through respectively $\bs_1, \bs_2, \bu$. The group multiplication $(G^{s_1, s_1}, \pist) \times (G^{s_2, s_2}, \pist) \to (G^{u,u}, \pist)$ is a Poisson morphism, and one can check that its restriction gives a Poisson isomorphism  
$(\Sigma_{\bs_1}, \pist) \times (\Sigma_{\bs_2}, \pist) \cong (\Sigma_\bu, \pist)$. One thus has
\begin{align*}
\Sigma_{\bv}&=\left\{\left(\begin{array}{ccc} p_1t_1 & -t_1 & 0 \\ t_1 & -q_1t_1 & 0 \\
0 & 0 & 1\end{array}\right)\left(\begin{array}{ccc} 1 & 0 & 0 \\
0 & p_2t_2 & -t_2 \\
0 & t_2 & -q_2t_2\end{array}\right): \; t_1^2(1-p_1q_1)=1, \, t_2^2(1-p_2q_2) = 1\right\}\\
& \cong \Sigma \times \Sigma =\{(p_1, \,q_1, \,t_1, \,p_2, \,q_2, \,t_2): \; (p_j, q_j, t_j) \in \Sigma, \;j = 1,2\},
\end{align*}
where $\Sigma$ is given in \eqref{Sigma-defn}, and $\pist$ is identified 
with the direct with Poisson bracket given in \eqref{bracket-Sigma}. 
On the other hand, parametrize
$BvB/B \subset G/B$ by 
$$
\Cset^2 \ni (z_1, z_2) \longmapsto [z_1, z_2] \,\stackrel{{\rm def}}{ = }\,
 \left(\begin{array}{ccc} z_1  & -1 & 0 \\ 1 & 0 & 0\\ 0 & 0 & 1\end{array}\right) \left(\begin{array}{ccc} 1  & 0 & 0 \\ 0 & z_2 & -1\\ 0 & 1 & 0\end{array}\right)_\cdot B \in BvB/B.
$$
The Poisson structure $\pi_1$ on $BvB/B$ is then given by $\{z_1, z_2\} = -z_1z_2$. One
checks that the groupoid structure on $\Sigma^\bv$ over $BuB/B$ is given as follows:
\begin{align*}
&\text{source map}: \; \theta(p_1, \,q_1, \,t_1, \,p_2, \,q_2, \, t_2)= [p_1, \; p_2t_1^{-1}],\\
&\text{target map}: \;\tau(p_1, \,q_1, \,t_1, \,p_2, \,q_2, \, t_2) = [p_1t_2^{-1}, \; p_2],\\
&\text{inverse map}: \; \iota(p_1, \,q_1, \,t_1, \,p_2, \,q_2, \, t_2) =(p_1t_2^{-1}, \; -q_1t_1^2t_2, \; t_1^{-1},
\; p_2t_1^{-1}, \; 
-q_2t_1t_2^2, \; t_2^{-1}),\\
&\text{identity bisection}: \; \epsilon(z_1, \, z_2) =(z_1, \,0, \,1, \,z_2,\, 0, \,1),
\end{align*}
and the groupoid multiplication is given by
\[
\mu(\gamma, \, \gamma^\prime) = 
 (p_1, \;\, q_1^\prime t_2^{-1} + q_1 (t_1^\prime)^{-2}, \;\, t_1t_1^\prime, \;\,  p_2^\prime, \;\, q_2^\prime + q_2  (t_1^\prime)^{-1}(t_2^\prime)^{-2}, \;\, t_2t_2^\prime),
\]
if $\gamma =(p_1, \,q_1, \,t_1, \,p_2, \,q_2, \, t_2)$ and $\gamma^\prime 
=(p_1^\prime, \,q_1^\prime, \,t_1^\prime, \,p_2^\prime, \,q_2^\prime, \, t_2^\prime)$ with
$p_1t_2^{-1} = p_1^\prime$ and $p_2 = p_2^\prime (t_1^\prime)^{-1}$. 
\hfill $\diamond$
\eex

\bre{re-concluding}
For any $v \in W$ and any representative $\bv$ of $v$ in $N_\sG(T)$, the symplectic groupoid $(\Sigma^\bv, \pist)$
over $(BvB/B, \pi_1)$ is algebraic in the sense that $(\Sigma^\bv, \pist)$ is an algebraic symplectic
variety and that the structure maps for the groupoid are all algebraic morphisms. However, 
as one can already see in the example of $SL(2, \Cset)$, the source fibers of these groupoids
are not necessarily connected. It would be interesting to understand how source-fiber connected symplectic groupoids
can be constructed from the ones in this paper.
\hfill $\diamond$.
\ere

\end{document}